\documentclass{notices}
\usepackage{amsfonts,amssymb,amsmath,amscd, xcolor, graphicx}

\title{Interview with Hyman Bass
}

\author{
  Hyman Bass
  \affil{
    Hyman Bass is a professor of the Marsal Family School of Education and the Samuel Eilenberg Distinguished University Professor of the  Department of Mathematics at the University of Michigan. His email address is hybass@umich.edu.
    }
  \and
 Lisa Carbone
  \affil{
    Lisa Carbone is a professor of mathematics at Rutgers University. Her email address is carbonel@math.rutgers.edu.
   }
    \and
Yvonne Lai
  \affil{
     Yvonne Lai is a professor of mathematics at the University of Nebraska-Lincoln. Her email address is yvonnexlai@unl.edu.
   }
}

\begin{document}

\maketitle

{\it 
Yvonne Lai$^\ddag$ and I$^\dagger$ had the privilege of interviewing my PhD advisor, Hyman Bass$^*$, about his remarkable and vibrant mathematical life, which has spanned seven decades. 

Hyman's mathematical research is in algebra and its relationship to algebraic geometry, number theory, topology, and geometric group theory. He was an inaugural Fellow of the American Mathematical Society in 2012. He received the U.~S.~National Medal of Science in July 2007 from President George W.~Bush. 

In the 1990s, Hyman began to investigate mathematical knowledge for teaching. His education research now also involves reasoning and proof in school mathematics, analysis of curriculum materials, characterizing mathematical and pedagogical practice across various contexts, and connection-oriented mathematical thinking.

Hyman has had 25 mathematics PhD students, with 174 mathematical descendants, and 3 mathematics education PhD students.

He has been elected to the National Academy of Sciences, the American Academy of Arts and Sciences, the National Academy of Education, and the Third World Academy of Sciences.  And he is a Fellow of the American Association for the Advancement of Science.

His mathematical writing is well-renowned for its systematic elegance. He won the Van Amringe Prize from Columbia University for the book {\it Algebraic K-theory} and the Frank Nelson Cole Prize in Algebra from the AMS for his work ``Unitary algebraic K-theory'' in Springer Lecture Notes in Mathematics, Volume 343, 1973. His paper ``Mathematics and Teaching,'' (Notices of the AMS, 2015) was selected for inclusion in {\it The Best Writing in Mathematics, 2016} published by Princeton University Press \cite{Pitici2017}.

He has a long and distinguished record of service to the profession. He was Vice President of the AMS in 1980-1981, President of the AMS in 2001-2003, and President of the International Commission on Mathematical Instruction in 1998-2006.  He was awarded the Mary P. Dolciani Award ``for distinguished contributions to the mathematics education of K-16 students'' by the MAA in 2013 and the Yueh-Gin Gung and Dr.~Charles Y.~Hu Award for Distinguished Service to Mathematics" by the MAA in 2006. He was Chair of the AMS Committee on Education from 1995 to 2000.}

\bigskip

{\it Acknowledgments. We are grateful to graduate student Em Stephen, who dedicated hours to wrangling a Zoom transcript into well-written natural language, with the support of the algebra faculty at the University of Nebraska-Lincoln especially Mark E.~Walker. We also thank T.Y. Lam for his careful review of the resulting manuscript, and Isaac Bass for reading this manuscript and providing photographs  with personal and family context.
}

\bigskip

\section{1930s-1959: Becoming a Mathematician}

\subsection*{In the Second World War}

{\it Interviewer: 
You were born 
with World War II on your doorstep. How did this shape you?}

\noindent {\bf Bass:}  
I was the seventh of eight children. My family was active, robust, and warm.

That was my world in early childhood.
The first set of conscious memories of being aware of the larger world came about the time of World War II.

\begin{figure}[ht]\begin{center}
		\includegraphics[scale=0.215]{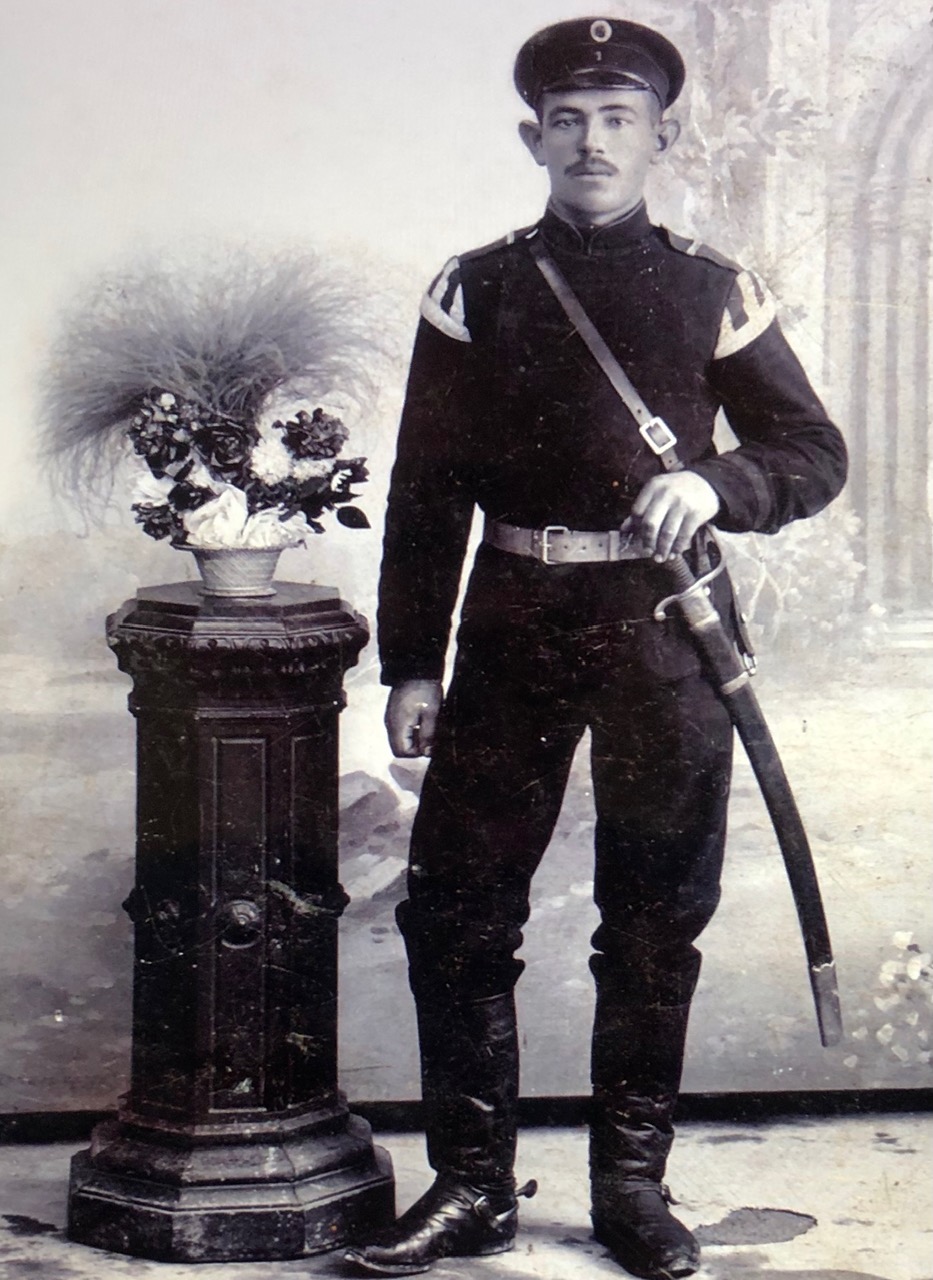}
        \caption{Hyman's father (Asriel Bashakevitz, anglicized in the USA to Isadore Bass) in the Russian Army, circa 1905,  Vilnius, Lithuania.
       }
		\end{center}
\end{figure}

My father, after emigrating from Lithuania (in 1911)\footnote{Hyman's mother Frume (anglicized in the USA to Fanny) Weissbord emigrated from Belarus in 1910.}, worked for an uncle who owned a kosher butcher shop in Houston.

  Later my father took that over and went into the wholesale meat
  business, which delivered to grocery stores and other markets around Houston.  

  \begin{figure}[ht]\begin{center}
		\includegraphics[scale=0.1645]{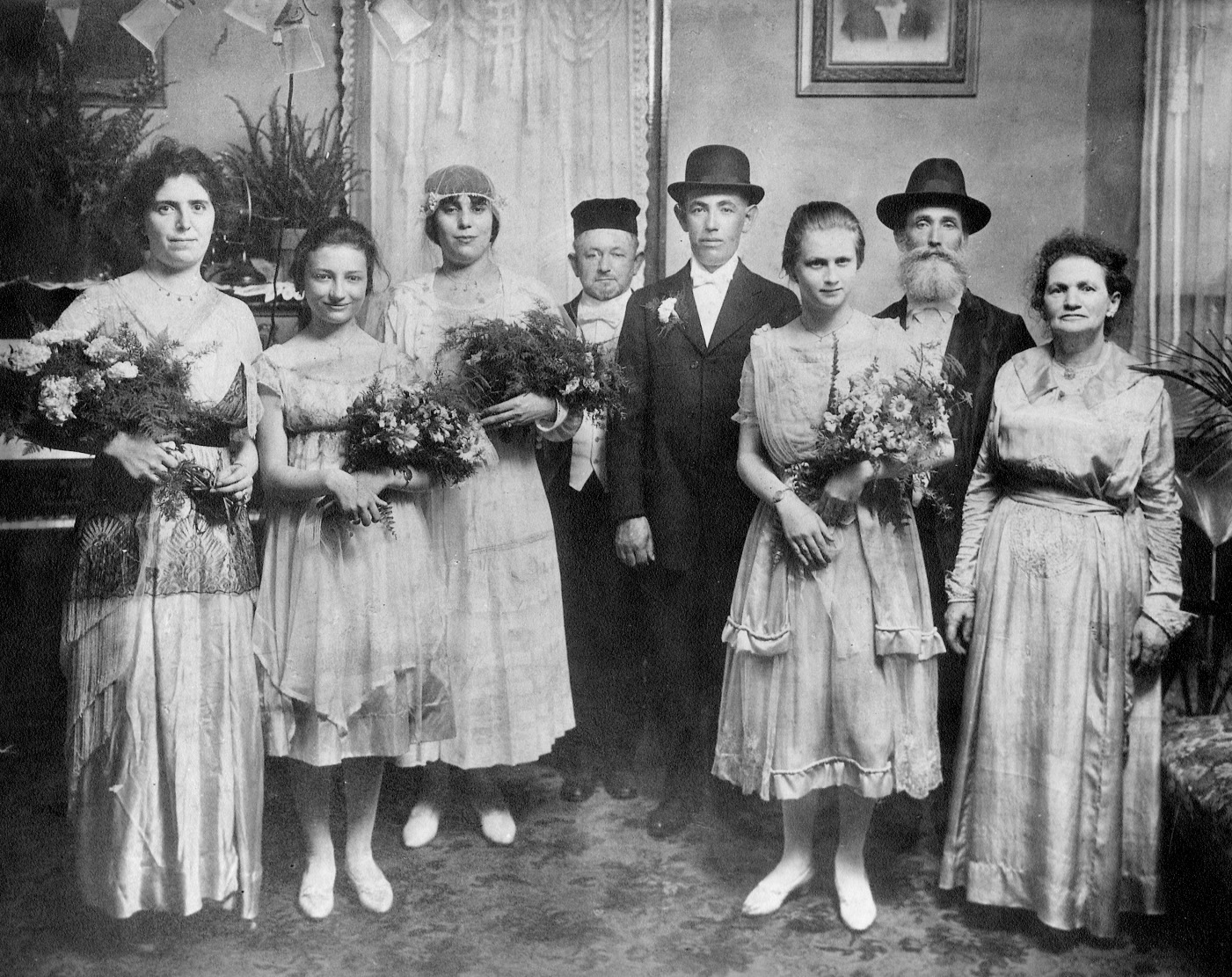}
        \caption{The wedding day of Hyman's parents Fanny Weissbord and Isadore Bass, 1917,  Houston, Texas.
       }
		\end{center}
\end{figure}

  I'm sure my father was looking forward to raising sons who would help him in what is a pretty physical business. 
  Then he and my mom rapidly produced four girls.  Finally, when the four sons arrived, the Second World War was beginning. 
  
  My two older brothers went to the Navy. My sisters went off to Washington. One joined the Women's Army Corps, and served in France shortly after the Normandy invasion.
  
 In grade school, my father had no help from my older brothers, so I would leave school about noon to help him on the delivery truck.  And so, my first experience working was carrying meat into markets to be weighed on scales. And
  my early education of the `real world' was
  delivering meat on a
  pickup truck, and working at a
  processing house
  where cattle were
  slaughtered, dressed, and put in coolers. 
  
  My companions on the loading dock were
  mostly African American men, who called me ``Little Bass.''
 I remember two brothers with whom I became friendly. I was always intrigued by their names: Napoleon Gibson, and Ulysses Gibson.
It took me a while to appreciate the significance of those names. 


 My brother Manuel -- whom we called Mangie in the family -- was in the Naval Officer Training Corps, which was basically engineering. Through this, he was learning a lot of science. This planted a seed of my own interest in science.
 
  Like all of my siblings, he was bright. He was also curious about the world, excited by learning, and eager to share everything he learned. 
When he came home on leave, he would teach my younger brother Isaac and me about the science he was learning. This planted a seed for an interest in science in general. I even did a project, possibly for Westinghouse, with  the fruit fly {\it drosophila melanogaster}. That was probably partly the influence of my brother again. But it didn't take hold the way mathematics eventually did.

 \subsection*{Becoming an algebraist in the post-war and Sputnik years}
 
{\it Interviewer: At what age did you know that you were very dedicated to mathematics?}

\noindent {\bf Bass:}  
   Not until college.

Near the time that my father retired from his business, we moved to California. Our summer vacations were often either in Colorado, in the mountains,  or else in
California. 
My dad was frugal, so in California, we were in Ocean Park, which is between Venice
and Santa Monica. We stayed in a small structure in the middle of an asphalt parking lot–but it was only a block from the ocean.

When my parents moved, one of my sisters, Silva, who was interested in the film industry, was already working there as a writer. She picked out the house that my parents bought. 
That's where I finished junior high and went to Hamilton High School, in LA.

At my high school, the natural course for people oriented towards science was to go to Caltech.
By then, Mangie was a graduate student at Caltech in geophysical science. My oldest brother Leon, who had been a radio operator in the Pacific, was there too, in electrical engineering. (My younger brother Isaac later went to UC Berkeley, earned a PhD in physics at Columbia, and now works on lasers at the Lawrence Livermore Lab.)

But I wanted more of a liberal arts education. When I asked Mangie for advice about this, he didn't have an immediate answer, but he consulted his classmates. They said Princeton had a good liberal arts program. I applied, and received something called a regional scholarship, perhaps for geographical diversity.

\begin{figure}[ht]\begin{center}
		\includegraphics[scale=3.28]{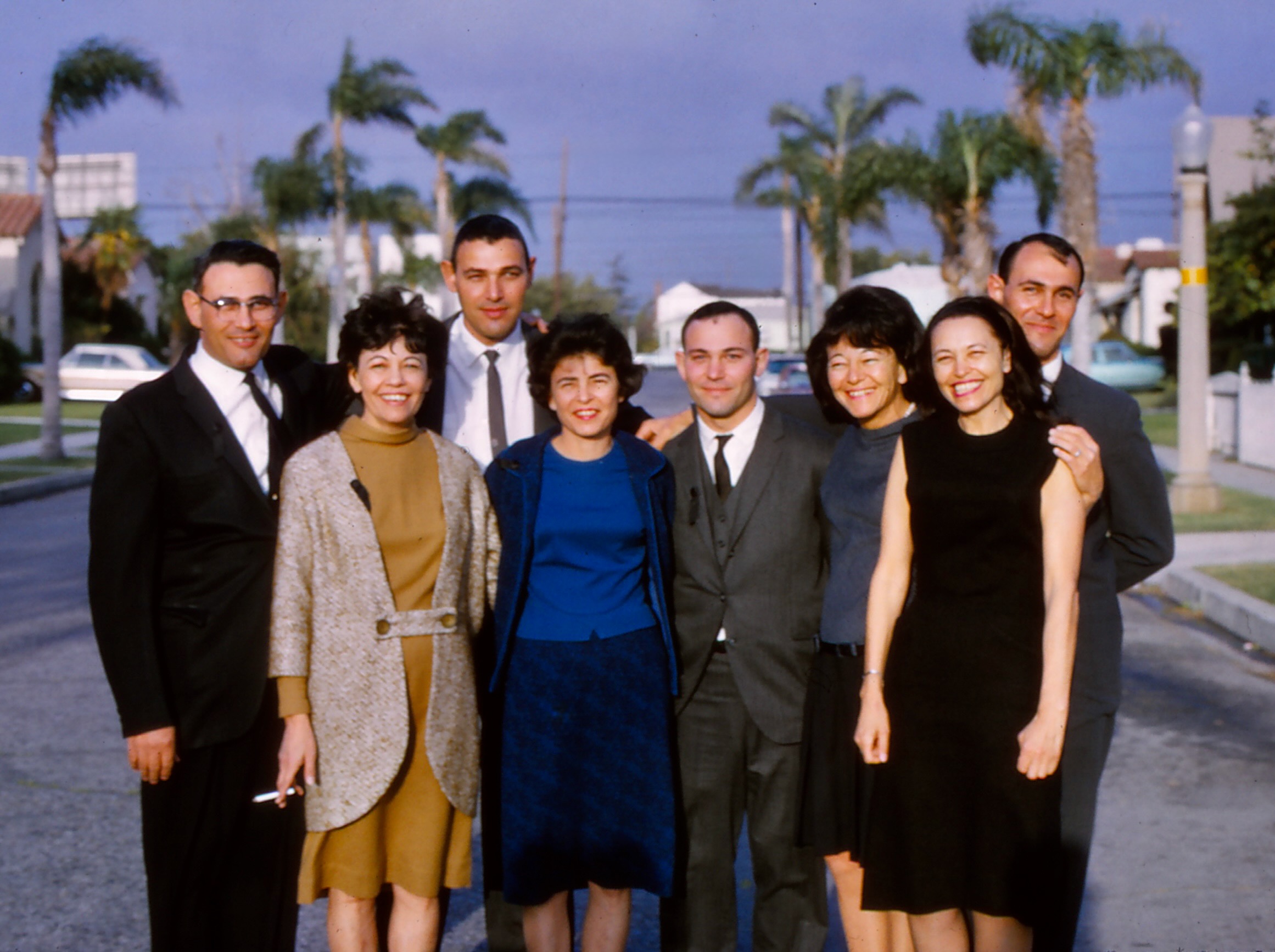}
        \caption{The Bass siblings left-to-right Leon, Pearl, Manuel, Frances,  Hyman, Sylva, Madeline and Isaac in Los Angeles, 1966.
       }
		\end{center}
\end{figure}

I didn't know anything about Princeton, and I inherited the frugality of my parents. When I went from California to Princeton, I hitchhiked. I made
several cross-country trips that way, not wanting to be a financial burden. The total expenditure my parents made on my education
was about \$500. I worked each summer. The summer between high school and Princeton, I worked in a logging camp in Oregon. It shut down because of fire danger, and then I worked as a forest fire fighter. The summer after freshman year, I worked on a salmon fishing boat near the San Juan Islands in Puget Sound.

When I arrived in Princeton, it was a shock. 
Coming from the California with blue skies and pastel colors, I found
earthy colors and leaded glass windows. I found it depressing at first, although eventually I came to like it.

Because I arrived a couple of days early, nobody was there. So I hitchhiked to Washington, DC. With no money to spend, I walked up and down the Washington Monument, then returned to Princeton. By
then, people were showing up.

You asked me when I turned toward mathematics. There's a very clear and definitive answer to that: freshman honors calculus, taught by Emil Artin. 
The
teaching assistants included John Tate and Serge Lang. The students in that class included Fred Almgren,
Mike Artin, Steve Chase, Lee Neuwirth, Steve Schanuel -- a number of people who became mathematicians. 

Emil Artin was a dramatic instructor with a stern countenance and operatic style.
On one occasion, a student said something in class. 
Artin raised his shoulders and walked toward the student. It seemed so menacing. But then he took a nickel out of his pocket, gave it to the student, and said, ``That was a good idea.'' 

At that time, calculus was not taught in high school. So
this was a first exposure for me, and 
it became my serious introduction to mathematics.

I always thought I knew what a real number was. Suddenly I was made to believe that it was an equivalence class of Cauchy sequences of rational numbers. That seemed like a real backward move. How is anyone going to work with an idea like that? 

It's not that I loved real analysis or calculus. I loved the kind of
mathematical thinking 
that I found not only in the instruction, but also among the students.
They would come up with these ideas, and I would think, {\it how on earth would anyone think of that?}
I guess if you're around that enough and you pay close attention to it, eventually you learn some of it, but not all.

\smallskip

\noindent {\it Interviewer: What did you major in?}

\noindent {\bf Bass:}  
Math, by default. I was interested in many subjects, 
but only in math did I feel competent.

\smallskip

\noindent {\it Interviewer: Did you have any other interests that
potentially challenged your interest in mathematics?}

\noindent {\bf Bass:}  I was curious about many things. But I never felt that any challenged my interest in mathematics, but actually in many ways enhanced it, especially my interest in art.

By the way, Princeton at the time was all male. 
And, it was expensive: if you had a date, you had put them up at a hotel. The best way I found to integrate into the social milieu was bartending at football weekends. Social performance consisted in
over-drinking alcohol. I didn't try to compete too heavily in that.

\smallskip
\noindent {\it Interviewer: When was it apparent that you wanted to go to graduate school?}

\noindent {\bf Bass:} That was taken for granted once I was a major. 
I loved 
mathematical thinking, and there was no
place else to do that, outside of graduate school.

Applied mathematics was not developed so
robustly as it is now.
The main application of mathematics
came during the Second World War. In fact, Samuel Eilenberg and Saunders MacLane began
their collaboration during
a war program at Columbia that was devoted to
ballistics among other issues.\footnote{This program was part of the Applied Mathematics Panel, created in 1942 as a division of the National Defense Research Committee within the Office of Scientific Research and Development \cite{Rees1988}.} 

Many of the people like von Neumann and others, who worked on the Manhattan Project development of the atomic bomb, were
mostly pure mathematicians, or physicists. 
Broad,
conceptual, and
big thinking
was important to the work.

Fields like applied math were of course important in the design of aircraft, and partial differential equations and dynamics were crucial to the design of high-tech equipment, especially some of the weapons in the war effort. These and logistics questions were eventually treated more amply in the post-war years in areas of applied mathematics.

One big effect of World War II was that because of the importance of radar and eventually atomic weapons and cryptography, the military
came to appreciate that fundamental research has a big role to play in national security.

They also understood that what was needed to make this happen, and they wanted it to continue after the war.
The idea was to set up some institutional way to continue to support that. There was a proposal for what became the National Science Foundation. But it was vetoed by Truman, who resisted its being highly independent of Federal oversight.

So the effort to support fundamental research continued, but in the Office of Naval Research.
When the NSF was finally established, it was modeled on the culture that existed then in the Office of Naval Research, which really supported basic science.

\begin{figure}[ht]\begin{center}
		\includegraphics[scale=0.74]{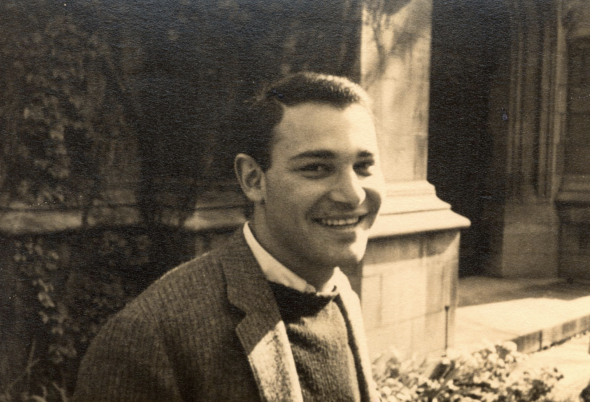}
        \caption{Hyman  as a PhD student at the University of Chicago, 1958, photo by Paul Halmos.}
		\end{center}
\end{figure}
From 1951 to 1955, I was in Princeton and then in graduate school in Chicago until 1959. In 1957, Sputnik occurred. 
A lot of people believed that
Russia was so backward that they would not be capable of anything like that. So it was a real awakening and shock.
The country 
discovered that they needed to take this seriously and to develop much more capacity and strength in technology and science. So they created incentives
to support
education moving in that direction. So any student
who intimated in any way that they were interested in science or math
was amply encouraged. So in my graduate and post-graduate years, it was easy to get a scholarship or fellowship, because of the support for progress in basic science and mathematics. 

\smallskip

\noindent {\it 
Interviewer: When did you realize that you were good at or particularly interested in algebra?
}

\noindent {\bf Bass:} 
From the beginning. Already in the calculus course, the thinking was algebraic.
We were not solving complicated differential equations; real analysis at the core is about metric spaces. 

There 
was one problem that I remember to this day.  It was from an exam in the course, which were all take-home exams where we had several days to work on them. Of course, there were no internet or modern calculators back then. 

\begin{quote}
{\it Define $f :\mathbb{R}\to\mathbb{R}$ as follows. When $x\in \mathbb{Q}$, and $x=\frac{p}{q}$ with $p, q\in \mathbb{Z}$ in lowest terms, we take $f(x) =\frac{1}{q}$. When $x \notin \mathbb{Q}$, we take $f(x)=0$. Where is $f$ continuous?} 
\end{quote}

My first reaction to this problem was anger. {\it How on earth can anyone answer a question like that?} You can't draw a graph. So I had to think, {\it what tools do I have to answer this question?} There were two tools: the definitions of continuity and of $f$. I wondered, {\it could you really reason just using those?}

It took
a couple of days to think through the problem. What I learned was not so much the answer, but the fact that
reasoning from definitions like that
was actually possible, and that you could answer
what seemed at first like an impossible question.

At first sight the function seemed like an artificial construct. But actually the problem embedded 
substantial ideas about discreteness, density,
and how they relate to each other. Those ideas
stayed in my mind and motivated many things I think about, like discrete groups and lattices.

\smallskip
\noindent {\it Interviewer: Why the University of Chicago for your PhD?}

\noindent {\bf Bass:} I forget exactly why I chose Chicago.

I should say, it's much harder now.
I don't think I would have succeeded as a student in today's world, where mathematical exposure seems so much more competitive and test-driven. 

I never once felt, at Princeton or Chicago, that I was competing with my classmates. 
The only competition was with the ideas. Mathematics was an adversary, but not a mean one, and the rules are clear. 
The triumph is not
putting something down, but rather
gaining
a mastery in the sense of a kind of companionship with the ideas.
That's one of the things about the culture of mathematics that I appreciated from the beginning.

Even in Princeton, 
there was
a whole range of personality types and physical types in that room.
The way they were treated only had to do with
the nature of their thinking.

My first published paper \cite{Bass58} came out of a U. Chicago course of Halmos on ``Algebraic Logic." Halmos had introduced the concept of a ``Monadic algebra." He mentioned that it was an open question whether a free monadic algebra on finitely many generators was finite or infinite. I showed that it was finite, and gave its cardinal. That was my first paper. It wasn't very hard.
\subsection*{Instruction and the language of mathematics in the sciences}

{\it Interviewer: If your younger self could time travel from your PhD years to the present day, what changes in the landscape would you be most surprised by?} 

\noindent {\bf Bass:} As a graduate student, I thought that mathematics had such a long and noble history, it would continue more or less in that way. I didn't foresee the radical, fundamental changes that have happened, 
influenced by the progress of the field (notably Grothendieck), and especially by technology.
And at various times, I have thought about the dynamics of the field as a whole, as a culture. 
Population size has made a big difference. 

When I was a student, American mathematics was beginning to burgeon. Before this, American mathematics had a slow start. Even someone like MacLane went to G\"ottingen for his degree. But then American mathematics flourished 
as it gained status
in science and math. 
And now, the number of people doing mathematics
has grown
enormously. Not everybody can become an academic. 
But, for a long time
there were plenty of academic positions open. The field had gained prominence, not only because people admired it, but from its growing impact on national security, on the economy, and on many other things. Even technology and applied areas depended hugely on mathematics.

Take something like
non-invasive medical diagnostics.
The big breakthroughs in that depended hugely on mathematics that was done at the end of the 19th century. 
It was on the shelf, ready to be used by chemists or medical practitioners.
For example, say you take a three-dimensional body, you take parallel slices, and you know data from those slices. 
Can you reconstruct the original three-dimensional body? 
That's a natural question in differential geometry.
Or if you take line sections from multiple directions, can you reconstruct it from those data? 
Those are geometric questions that were largely solved by the late 1800s. 
They are the basis for
CAT scans and
MRIs. 

There is a beneficial paradox for mathematics, that Eugene Wigner called the ``unreasonable effectiveness of mathematics'' (1960, p. 1 \cite{Wigner1960}). Over and over again, 
ideas of mathematics  are pursued for
reasons of intellectual curiosity.
But it turns out that our curiosity and aesthetic instincts
lead us to the same answers that nature has found for many questions. 
It may be a long kind of arc before these things retouch, but the point is when they do,
they produce ideas that someone who was doggedly focused only on the application would never have arrived at.
They would never think
of the
stages of abstraction and conceptual development
that mathematics evolved to consider. 

As Galileo said, we have to understand the language of nature, and that is the language of mathematics.\footnote{
``Philosophy is written in this grand book, the universe \dots [which] cannot be understood unless one first
learns to comprehend the language and read the letters in
which it is composed. It is written in the language of mathematics\dots '' Galileo/Drake, 1623/1957, pp.~237-238\cite{Galileo1623}
} 
It is a deep insight, and it keeps getting reproven. It is a fortuitous 
repetition of history that the
aesthetics of the subject keeps touching ground with the real world. 

I think we've kind of saturated the academic
pathways into mathematics.
More and more people with mathematical training have to take jobs in finance, data science, computer science, and
other applied areas.
That's fine, because
good critical thinkers are important in those areas.
The demands of the so-called real world are increasingly in need of people who are
good thinkers and problem solvers; mathematics can develop these skills. 

\smallskip

\noindent 
{\it 
Interviewer: What do you think has improved for mathematics since you were a postdoc? And what do you think has become more challenging?}

\noindent {\bf Bass:} Generally the field continues to progress, and I don't think that will ever cease
as long as we don't destroy civilization, which is unfortunately something we can't rule out anymore.

Here's an example. What's the economic foundation of academic mathematics and of the sciences in general?

In the sciences, it's running labs. If you are a department chair in the sciences, and you appoint a new faculty member, you set them up, you declare a salary, and in addition there has to be what used to be half a million dollars, and maybe much more now, just to set up their laboratory. Labs not only need academic support. They need staff, management, and administrative skills, even from the scientists, to keep projects running. And, accordingly, research grants in science are, in general, much larger than those in mathematics.

On the other hand, mathematicians have a different economic foundation. I think that for mathematicians, the lab is the classroom -- especially advanced graduate courses that act as seedbeds for ideas in the subject. The culture of mathematics has a different structure than it does in science. In an active mathematics department, the list of seminars has no counterpart in other fields. Huge time and intellectual effort go into these. Faculty don't typically receive formal credit for running seminars. 

Mathematics departments typically have one of the larger faculties in the STEM fields.
What allows mathematics faculty 
to have the size that it does? 
The economic foundation is educational: the instructional mission.
Mathematics is an enabling discipline for sciences and technology. When foundational subjects such as calculus, linear algebra, and differential equations are taught, mathematics has jurisdiction over the instructional terrain.
And as universities depend more and more
on tuition for income, instruction becomes an increasingly crucial resource. 
It sustains reasonable faculty size, and it's very precious. 

Mathematicians do not sufficiently appreciate
how much of a stake they have in the quality of what they do
with undergraduate instruction.
It's basic economics. Other people, especially the engineers, often covet that territory,
and rightly so.
They feel that they are far better equipped to teach
calculus to engineers, because they know how it will be used. This argument is not without merit. 

Eilenberg was chair of the department at Columbia when such a challenge occurred. 
Bypassing normal process, he did not go to the Dean to argue the case. Instead, he went straight to the University President. 
He said that Columbia is a first rate research university, and a research university without a first class mathematics department is not tenable. A first class mathematics department needs a reasonable faculty size that can cover enough areas with deep mathematical strength.

Now, Columbia does not have a huge math department, but for its size, it has relatively a large instructional role compared to the size of the university and the college. And if they lose that instructional real estate, it means a critical reduction of the department size, and therefore  its quality. It is a question of critical mass. This argument
was persuasive.

But this argument is not purely about the merit of the mathematics. It also carries a message
that
the quality of the instruction has to justify the size. Unfortunately, this is not always sufficiently communicated to the young mathematicians who come into the
department and don't realize how much they have at stake in it and in the quality of the instruction.

This has changed slowly and in good ways over time. You can no longer say that if someone is doing really outstanding mathematics, it doesn't matter how inept they are as a teacher. The argument can further be made that, in fact, improving the quality of instruction will improve the quality and productivity of the mathematics.  This is a commitment that must be a part of department culture, and must be nurtured and supported by department leadership.

\section{1960s-1990s: Mostly commutative homological algebra, projective modules, ubiquity, and trees}

{\it Interviewer: In the 1960s to 1970s, your work was mostly focused on homological algebra and projective modules. What stands out to you in your memory about the the work you did in those years?}

\noindent {\bf Bass:}  
In preparation for this interview, I was looking through the book from my 65th birthday conference \cite{LamMagid1997}. I remember on a couple of occasions I read about some results I admired, and I asked  somebody who had proved that, and they said, ``You did." That's the way I feel about reading through that volume. A lot of the things have faded in my memory, except for very general kinds of impressions.

I do remember that it's always nice to come into an area that's fresh,
which homological algebra was at the time. The territory hasn't been completely tilled, and all the low hanging fruit hasn't been picked. It offers many new questions. Homological algebra was modeled on things that were developed in in topology, the homology and cohomology theories.

Saunders MacLane had this paper about natural transformations which were basically morphisms between functors. For example, a finite-dimensional vector space is isomorphic to its dual. They have the same dimension. But it's isomorphic to its second dual in a natural way, in that you can identify it with its second dual. Among those two isomorphisms, one is canonical; it's {\it natural}. And so, the notions of category theory began to emerge.  

Homological algebra, which contains the seeds of category theory, was just beginning to be developed as a novel theory that axiomatized ideas of algebraic topology. Topology is about complicated, continuous objects. Discrete things are easier to manage. They are more susceptible to calculation, making it easier to detect differences. For example, (co)homology theories attach discrete algebraic objects to spaces. So these theories were powerful tools, because they could algebraicize complicated geometric objects. Homological algebra axiomatized the methods of these theories. And axioms are then open to application to other categories than topological spaces, like rings, groups, Lie algebras, etc.

Descartes built a major bridge between algebra and geometry.  This has evolved to the practice of associating a space $X$ with a suitable ring $R(X)$ of functions on $X$.  This supports a partial dictionary for translating many geometric objects related to $X$ into algebraic analogues related to $R(X)$.  For example, a vector bundle on $X$ corresponds to a projective $R(X)$-module.  And projective modules have purely algebraic meaning for any ring, not necessarily of the form $R(X)$.  Thus, projective modules, and the homological algebra of rings are purely algebraic theories that have important links to geometry, but also connections, some unanticipated, to classical questions in algebra.

The new homological algebra of rings opened its door for application.  The first homological invariant of a ring $R$ is its global dimension.  I remember that Eilenberg investigated this for the most highly developed part of ring theory, that was, at the time, finite dimensional algebras. They have a nilpotent radical, modulo which they are semi-simple. It turned out that the global dimension of $R$ is zero if $R$ is semi-simple, and infinite otherwise.  So, disappointingly, homological algebra illuminated nothing new.
 
But a more productive site of application was in commutative algebra, which is tied to algebraic geometry. What does homological algebra have to say there? Here, there were huge payoffs. These were the developments of Jean-Pierre Serre, Maurice Auslander, David Buchsbaum, and other algebraists at the time. For example, finiteness of global dimension was equivalent to regularity. Geometrically, it was equivalent to smoothness of the corresponding variety. That was the birth of commutative homological algebra, which remains a vibrant area of research to this day. 

Irving Kaplansky was beginning to expose homological methods in commutative algebra in a series of lectures that inspired the direction of my dissertation.\footnote{Irving Kaplansky was Hyman's doctoral advisor.}  For my purposes, Kaplansky was a very comfortable kind of teacher. He did not try to cover the entire waterfront. He would pick selected key theorems, and then he would try to follow a geodesic to those theorems and not allow technicalities and side issues to get in the way. It turned out to be an appealing way of teaching. It sketched big ideas of the subject, and also laid open many unanswered questions at the time. For instance, consider the question of whether a regular local ring – the local ring at a smooth point – is a unique factorization domain. That's a perfectly natural algebraic question whose answer was proven using homological methods. Serre framed the theory of intersection multiplicities in homological terms. 

In its early stages, homological algebra brought into the foreground things about projective and injective resolutions, and projective modules in particular. Projective modules are very nice. Free modules are like vector spaces, but with scalars from a ring instead of a field. So then the idea is that you'd like to imitate linear algebra or elementary divisor theory for principal ideal domains, and see to what extent you can generalize that.

In my thesis I was investigating things to do with injective dimension. That was less actively treated at the time, so it was a little-worked territory for me. One reason, I think, that it was not given as much attention is that injective modules tend not to be finitely generated. For instance, for abelian groups, an important injective module is $\mathbb{Q}/\mathbb{Z}$, which is definitely not finitely generated.  So, the development of injective modules and injective dimension led me in a direction that seemed like a very natural thing to look at homologically. At the same time, it was not so clear how interesting it was more broadly in mathematics. I didn't worry about that at the time. But my so called ``Ubiquity'' paper, that culminated my work on this topic, ended up as one of my most cited papers  \cite{Bass1963}.

Given my approach, it was natural to investigate rings $A$ such that $A$ itself as an $A$-module has finite injective dimension.  I found several conditions equivalent to this, especially when $A$ is commutative Noetherian.  When I showed this to Serre, he pointed out that, when $A$ is an affine ring as in algebraic geometry, some of these conditions, related to duality, had been studied by Grothendieck, who called them Gorenstein rings, in reference to Gorenstein's thesis about plane curve singularities.  Other kinds of examples were integral group rings of finite abelian groups, relating to my interest in integral representations of finite groups.

The ubiquity paper was an effort to coherently assemble a host of not obviously related results, before turning my attention to other things.  This was a kind of ``house cleaning''; my aim was not to solve a big problem. And so, I found it curious that this paper drew so much attention.
 
To bring things up to date, we can look at group actions on trees.  This area emerged from work on the congruence subgroup problem, which asks, for an arithmetic group over a global field, to what extent are finite index subgroups ``explained'' by congruence subgroups.  In a split simple algebraic group $G$, effective results were found when $\textnormal{rank}(G) \geq 2$.  It was long known that this fails for rank 1, for example $SL_2(\mathbb{Z})$.

Serre went on to analyze the case of $SL_2(A)$, where $A$ is an $S$-arithmetic ring, for example $\mathbb{Z}[\frac{1}{p}]$, $p$ prime.  He showed that the congruence subgroup property is essentially rescued if the unit group of A is infinite.  His methods involved the study of discrete subgroups $\Gamma$ of $SL_2(F)$, where $F$ is a local field. The one notable result about this was Ihara's Theorem: {\it If $\Gamma$ is torsion free, then $\Gamma$ is a free group.} Serre systematized this study by using the action of $\Gamma$ on the Bruhat-Tits building, which, in rank 1, is a tree.  

In 1968, Serre gave a Coll\'ege de France course on these ideas. I was on a year's leave in Paris, and Serre invited me to prepare notes of his lectures. I was happy to do this, though my focus at the time was on algebraic $K$-theory.  Serre initiated a general study of a group $\Gamma$ acting on a tree $X$.  Beyond recovering Ihara's Theorem, he showed, for example, that if $\Gamma\backslash X$ is an edge, then $\Gamma$ admits an amalgamated free product decomposition.  In a caf\'{e} discussion about the course notes, Serre indicated how these results should be fragments of a general combinatorial theory of group actions on trees, and he proposed that I might work out the details in the notes. I did this, formulating it as a kind of analogue of ramified covering space theory, in which $\Gamma$ is represented as the ``fundamental group'' of the quotient ``graph of groups''  $\Gamma\backslash\backslash  X$. This work, now called ``Bass-Serre theory,'' was incorporated into Serre's book, {\it Trees}, based on the course \cite{SerreBass1977}\cite{Serre1980}\cite{Bass1993}. 

This done, my attention returned to algebraic $K$-theory, where, at the Institut des Hautes \'{E}tudes Scientifiques, I did some work with Tate on $K_2$ of global fields.
About a decade later, Peter Shalen asked me a specific question about finitely generated subgroups of $GL_2(\mathbb{C})$ that he suspected could be answered using the theory of group actions on trees.  I proved a little theorem confirming his suspicion.  When Thurston came to Columbia for a colloquium talk. I showed him the theorem, and asked if he knew what Shalen's interest might be.  He said that he didn't know, but that the result could help to prove the Smith Conjecture, about the un-knottedness of the fixed circle of a periodic diffeomorphism of the 3-sphere. Thurston's impulse turned out to be true, but its realization had to mobilize the work of several mathematicians, eventually assembled in the book, {\it The Smith Conjecture} \cite{MorganBass1984}.

This work led Morgan and Shalen to deepen and generalize the theory of group actions on trees.  In particular, they introduced the notion of a $\Lambda$-tree, where $\Lambda$ is any totally ordered abelian group.  Combinatorial trees are the case $\Lambda= \mathbb{Z}$.  They established the geometric importance of $\mathbb{R}$-trees \cite{MorganShalen1984}. These events prompted my reawakened interest in tree actions, including, in collaboration with Roger Alperin, the case when $\Lambda$ is non-archimedian. One outcome of this work was my book with Alex Lubotzky, {\it Tree Lattices}, written with the assistance of my students, Lisa Carbone and Gabriel Rosenberg \cite{BassLubotzky2001}. This book was a product of a long and fruitful collaboration with Lubotzky on geometric methods in group theory.  Alex's presence always seemed to energize and animate one's thinking. The work from this time such as by Carbone, Kulkarni, Rosenberg 
pretty much saturated what you could say about tree lattices and how they were analogous to the case of Lie groups \cite{BassCarboneRosenberg2001}\cite{CarboneRosenberg2003}\cite{Carbone2001}\cite{BassKulkarni1990}.

\begin{figure}[ht]\begin{center}
		\includegraphics[scale=0.50]{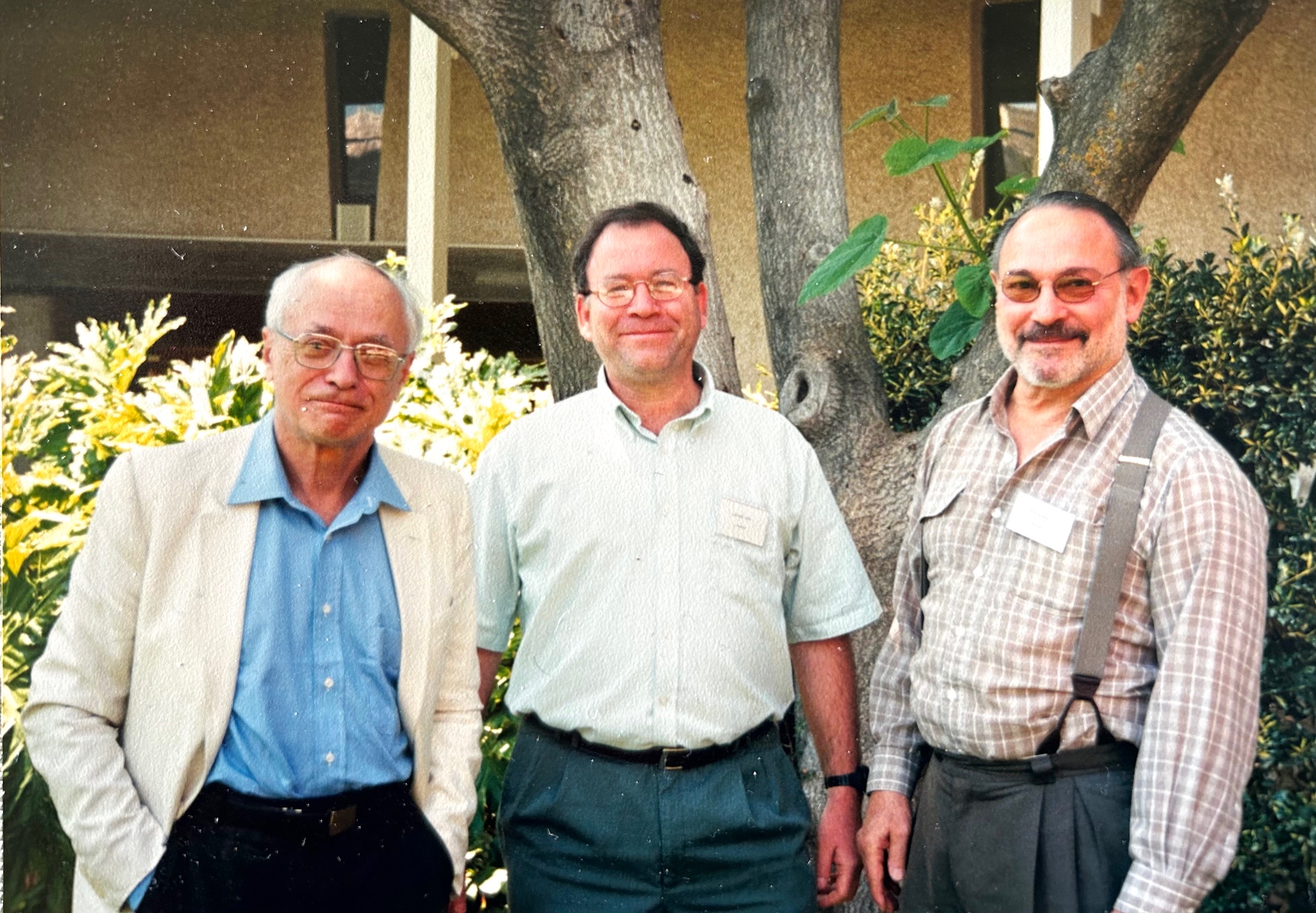}
        \caption{Left-To-Right J.-P. Serre, Alex Lubotzky and Hyman Bass, Einstein Institute of Mathematics, Hebrew University of Jerusalem, 1998, photo by Lisa Carbone.}
		\end{center}
\end{figure}

Ihara's paper, a key source of all this work, not only showed that a torsion free discrete subgroup $\Gamma$ of $SL_2(\mathbb{Q}_p)$ is free, but it attached to $\Gamma$ a zeta function \cite{Ihara1966}. Though I had not worked much with zeta functions, I decided, for my own edification, to understand that part of Ihara's paper. I generalized Ihara's construction to any uniform lattice, not necessarily torsion free,  on a not necessarily regular tree, and used ideas about non-commutative determinants arising from work on projective modules and representation theory, as, for example, in my paper on ``Euler Characteristics and characters of discrete groups'' \cite{Bass1976}. Though I treated this paper as an exercise in self-instruction, it became, again surprisingly, one of my most cited papers.

\section{Mid 1960s-mid 1980s: Algebraic K-theory and the congruence subgroup problem}

{\it Interviewer: It may be that what you are
 best known for is work on algebraic K-theory and the congruence subgroup problem from the 1960's to the 1980's. What do you recollect about this work? }
 
\noindent {\bf Bass:} 
It's interesting how things have a continuity to them.  
I had been thinking a lot about Serre's Conjecture, that projective modules over a polynomial ring are free. 
The model that I always had in mind as an algebraist was elementary divisor theory, that is, modules over a principal ideal domain or over a Dedekind domain. 

The big awakening was Serre's paper in one of the Dubreil-Pisot Seminars, where he proved that a projective module of sufficiently large rank has a free direct summand \cite{Serre1958}. This result was unlike anything that I had seen before. So I studied that proof very closely.

His proof made me realize that there was a big wall between two kinds of methods, topological and algebraic. The wall is akin to that between smoothness and analyticity. With analyticity, if you know a germ of a function – if you know it in a small neighborhood of a point – it is determined globally. Two analytic functions that agree in a small open set are the same. If it's only differentiable, even infinitely differentiable, that's no longer true.

This means that with smoothness -- and not analyticity -- you can make partitions of unity. You can make pieces of a function that behave nicely in different places, and you can splice them together smoothly.

Partitions of unity were the bread and butter of topological methods. In terms of viewing projective modules like vector bundles, getting a free direct summand means you can get a section which is non-zero everywhere, because then it will generate a 1-dimensional direct summand.

Serre's result is the algebraic analogue of a simple topological theorem: {\it If the rank of a vector bundle $E$ exceeds the dimension of the base, then $E$ has a (continuous) non-vanishing section.} So how do you make a section which is non-zero everywhere? You build it up inductively on the skeletons of the base space. Inductively, having a non-vanishing section on the
$t$-skeleton, you want to extend it to the interior of a cell of the $(t+1)$-skeleton. It has been constructed on the boundary, and you want to extend it without zero to the interior. If the bundle rank is $r$ and the dimension of the base space is $d$, the obstruction to this extension lies – almost tautologically – in $\pi_t(S^{r-1})$, which vanishes if $r-1 > d-1 \geq t$.

Functions in algebra behave like analytic functions. If you have two polynomials, even in many variables, and they agree on an open set, then they're the same. So you can't stitch polynomial functions given locally. 

In a projective module $P$ of rank $r$, Serre sought to build an element (``section'') $s$ that is non-zero modulo every maximal ideal. So one is working with the space $X$ of maximal ideals of the ring. The zero locus $Z(s)$ of an $s$ is closed in $X$. Choose a section $t$ non-zero on $Z(s)$ and ‘linearly independent' of all the sections used to build $s$. This is possible if $r$ is ‘large,' and then $Z(s + t) \subset Z(s)$, and has lower dimension. That is Serre's method. This was a totally different kind of argument than what I had seen before.  

Let $\mathcal{P}_r(R)$ be the set of isomorphism classes of rank $r$ projective $R$-modules, and define 
$\mathcal{P}_r(R)\longrightarrow \mathcal{P}_{r+1}(R)$ by sending $P$ to $P\oplus R$. With $d = \textnormal{dim}(R)$, Serre's theorem says that these maps are surjective for $r \geq d$. Inspired by the topological analogue, I was able to show that they are injective for $r > d$. This brought increased attention to ``stable range'' phenomena in algebra.

Then came Grothendieck's Riemann-Roch Theorem, fortunately published in a seminal paper of Borel and Serre \cite{BorelSerre1958}.\footnote{And later, it was published in SGA 6 \cite{Grothendieck1957}.} Grothendieck wanted to wait until the ideas were further developed, but the ideas were too fertile to delay publication.

One of these ideas is the $K$-functor, the ``Grothendieck group'' $K(X)$ of vector bundles on an algebraic variety or scheme $X$, treated as something like a cohomology theory. Atiyah and Hirzebruch ran with this and launched topological $K$-theory for topological spaces $X$, making 
$K^n(X)$ equal to ``$K$ of the $n^{\textnormal{th}}$ suspension of $X$''.

\begin{figure}[ht]\begin{center}
		\includegraphics[scale=0.066]{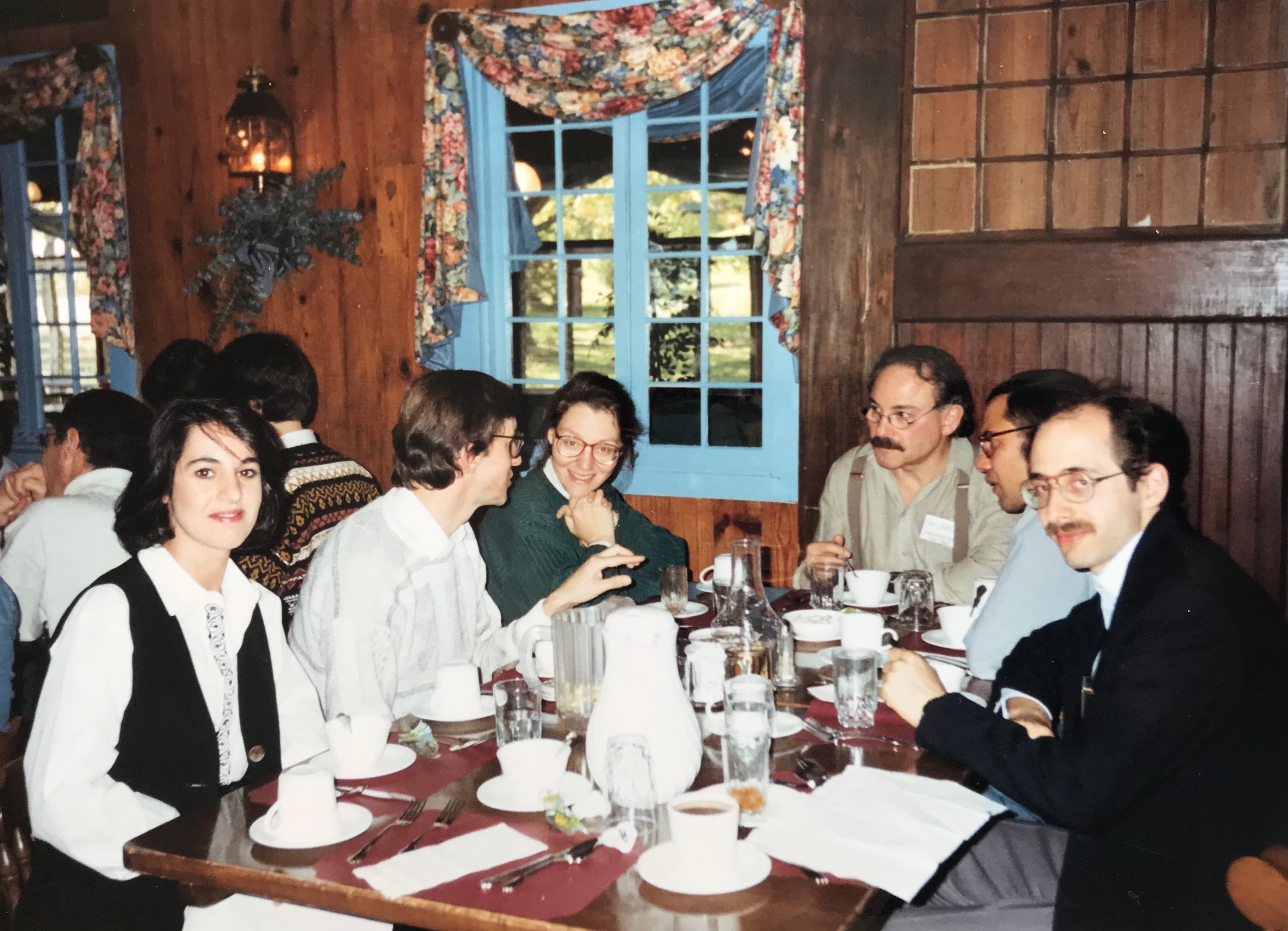}
        \caption{Left-to-right, Lisa Carbone, Peter Kropholler, Elizabeth Schneider, Hyman Bass, Sal Liriano and Ilya Kapovich, Albany Group Theory Conference, 1993.}
		\end{center}
\end{figure}

To me it then seemed natural to try to construct some algebraic analogue of these developments. For a ring $R$, one starts with $K_0(R)$, the Grothendieck group of finitely generated (left) $R$-modules. When $R$ is commutative, Noetherian, and regular, Grothendieck proved a kind of ``homotopy invariance:'' $K_0(R) \longrightarrow K_0(R[t])$ is an isomorphism. Thus, a finitely generated projective $R[t]$-module is ``stably isomorphic'' to $P[t]$, where $P$ is an $R$-module. Combined with the stability theorems, this gave partial affirmation of the ``Serre Conjecture'' that, when $R$ is a field, projective 
$R[t_1, t_2, \dots,  t_d]$-modules $P$ are free; this now followed from homotopy invariance and stability theorems if $\textnormal{rank}(P) \geq d + 2$. The Serre conjecture was fully proved about a decade later by Quillen and Suslin.

\begin{figure}[ht]\begin{center}
		\includegraphics[scale=0.092]{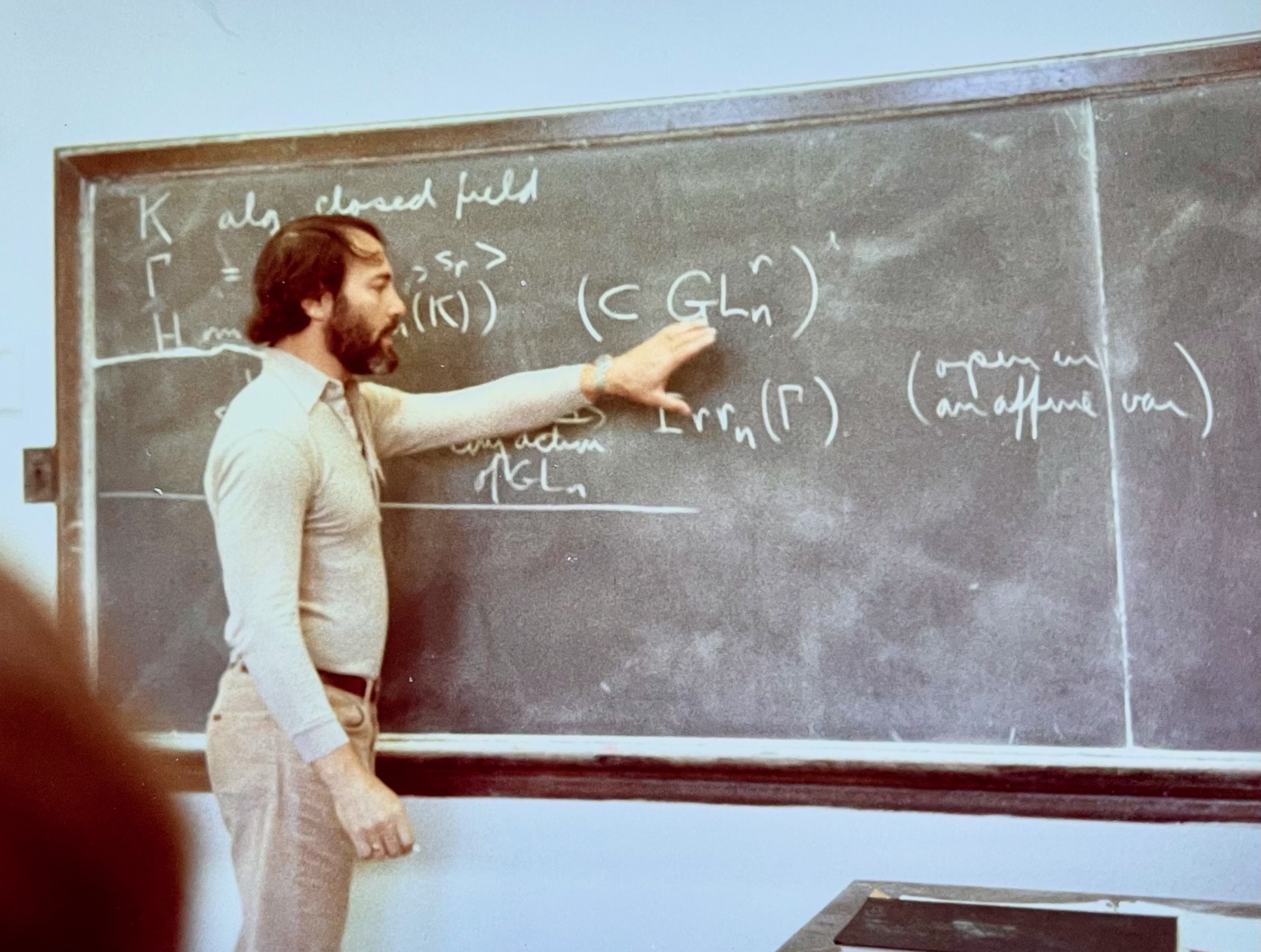}
        \caption{Hyman, mid-to-late 1960's}
		\end{center}
\end{figure}

It was natural to seek a construction of higher algebraic $K$-groups, $K_n(R)$, $n \geq 0$, but there was no satisfactory algebraic analogue of the suspension. I looked directly at the case $n = 1$. Topologically, $K_1(X) = K(SX)$, the Grothendieck group of vector bundles on the suspension $SX$ of $X$. $SX$ can be viewed as the union of two cones over $X$, glued along their common base $X$. Let $E$ be a rank $n$ vector bundle on $SX$.  Since the cones are contractible, $E$ can be identified with the trivial rank n bundle on each cone. Then $E$ is defined by a ``gluing'' automorphism of the trivial bundle on the ``equator'' $X$ of $SX$, in other words, by an element $ \alpha \in GL_n(C(X))$, where $C(X)$ is the ring of continuous real functions on $X$. Up to isomorphism, $E$ is determined by the homotopy class of $\alpha$, i.e., the image of  $\alpha$ modulo the identity component $GL_n^0(C(X))$ of $GL_n(C(X))$, so $K_1(X) = \lim_n GL_n(C(X))/GL_n^0(C(X)) = GL(C(X))/GL^0(C(X))$.

To define an algebraic analogue of this, one needs an algebraic proxy for $GL_n^0(C(X))$, with $C(X)$ replaced by any ring $R$. For various reasons it was found reasonable to use the group $E_n(R)$ generated by the elementary matrices $I + te_{ij}$, with $t \in R, i \neq  j$ and $e_{ij}$, $ (1 \leq i, j \leq n)$, is the standard basis of the space of $n \times n$ matrices. Further, define $K_1(R) = GL(R)/E(R)$. I proved some stability theorems for $K_1$ analogous to those for $K_0$. With Alex Heller and Dick Swan, we proved the analogue of homotopy invariance for $K_1$. In my book, Algebraic $K$-theory, I constructed some exact sequences relating $K_0$ and $K_1$ \cite{Bass1968}.

Further, Milnor wrote a beautiful introduction to algebraic $K$-theory in which he extended the theory to $K_2(R)$ \cite{Milnor1969}.

Though the rings $C(X)$ mediated the connection of the algebra to topological $K$-theory, it turned out that an entirely different part of topology had, much earlier, developed an interest in a slight variant of $K_1(R)$, when $R$ is the possibly non-commutative integral group ring $Z\pi$ of a group $\pi$. This was ``simple homotopy theory,'' asking whether a homotopy equivalence $X\longrightarrow Y$ is ``simple.'' J.~H.~C.~Whitehead defined an obstruction that lives in what we now call the Whitehead group $Wh(\pi)$, where $\pi$ is the common fundamental group of $X$ and $Y$. 
We can write $Wh(\pi) = K_1(Z\pi)/[\pm \pi]$, where $[\pm \pi]$ denotes the classes of diagonal matrices with entries in $\pm \pi$. Whitehead established the basic algebra of elementary matrices, showing for example that $E(R)$ is the commutator subgroup of $GL(R)$, and his student, Graham Higman did the first substantial calculations of $Wh(\pi)$, for finite abelian $\pi$. This algebra provided the basic tools for $K_1$ calculations.

For $R$ commutative, the determinant decomposes $K_1(R) = R^{\times} \oplus SK_1(R)$, $R^{\times}$ being the unit group of $R$ and the latter coming from the special linear group. Moreover there are relative groups $SK_1(R,J)$ defined for ideals $J$ of $R$. And already for $R = \mathbb{Z}$, calculation of the groups $SK_1(Z, J)$ encountered the century old ``congruence subgroup problem'' for $SL_n(\mathbb{Z})$ with $n \geq 3$. This enlisted the efforts of several mathematicians, notably Milnor and Serre, that culminated in our paper on the Congruence Subgroup Theorem \cite{BassMilnorSerre1967}. This resolved the $SL_n$ congruence subgroup problem $(n \geq 3)$ for arithmetic rings, and related it to reciprocity laws in number theory.

These diverse connections (to algebra, algebraic geometry, topological $K$-theory, simple homotopy theory, functional analysis, number theory) showed that these still modest developments of algebraic $K$-theory had more mathematical promise than I had anticipated. It intensified the interest in constructing a ``good'' higher algebraic $K$-theory. So it seemed like algebraic $K$-theory had lots of intriguing tentacles, some extending beyond my expertise, but no mathematical center. This prompted me to organize a conference assembling all of the diverse clients of this emerging subject. In the proposal to the National Science Foundation for this conference I said something like, ``There are lots of people with these related interests, but they're in very different areas. Just bring them together and let the human chemistry work.'' It was a spectacular success, thanks to Quillen, who brought with him a fully developed construction of higher algebraic $K$-theory (two versions), together with fundamental computational tools. This won him the Fields Medal. The conference proceedings were published in three thick volumes of the Springer Lecture Notes in Mathematics \cite{Bass1973one}\cite{Bass1973two}\cite{Bass1973three}.

\section{1970s-1980s: Writing with Bourbaki}

\begin{figure}[ht]\begin{center}
	
		\includegraphics[scale=3.20]{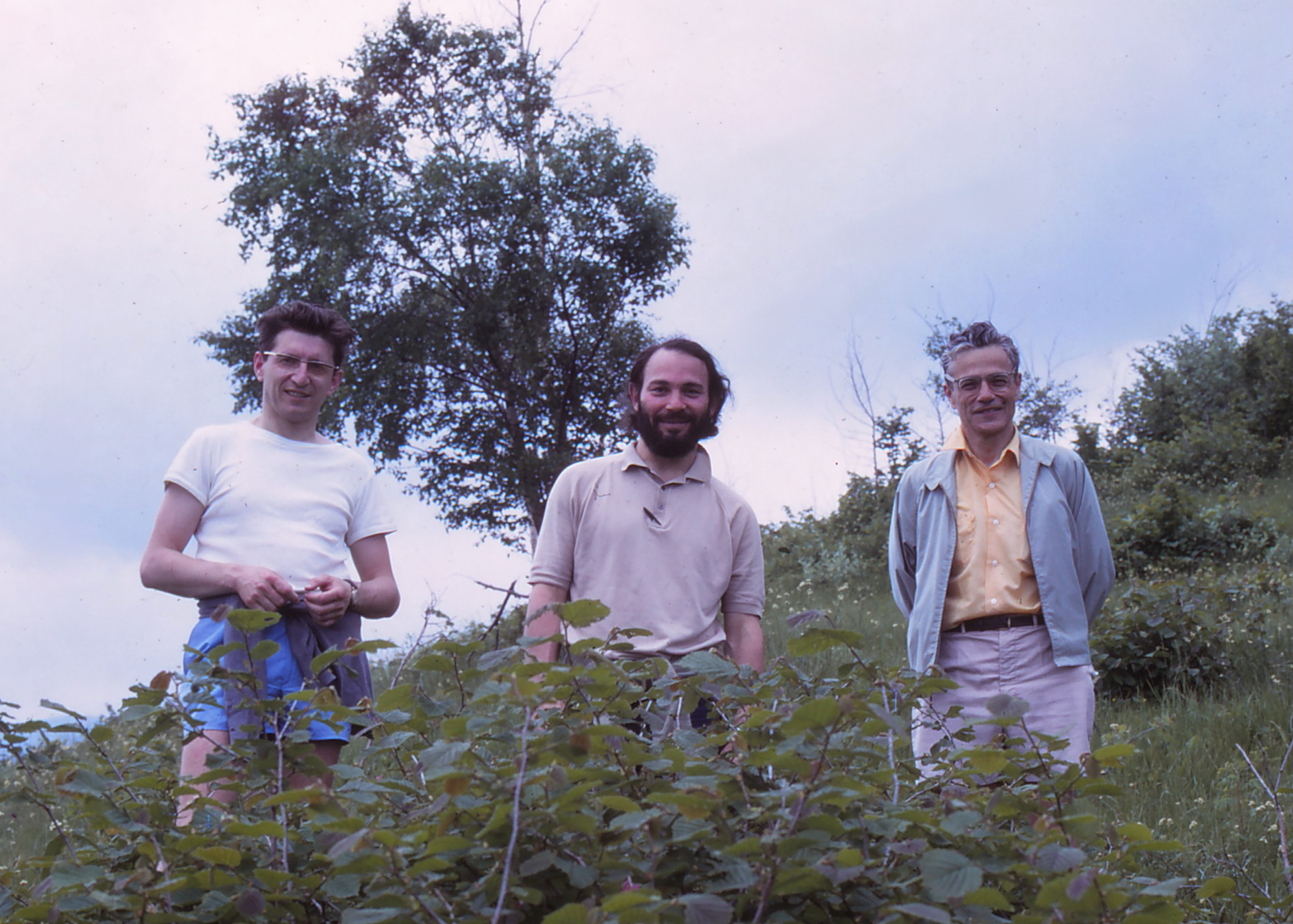}
        \caption{Left-to-right Michel Raynaud, Hyman Bass and Armand Borel, Corsica, 1970's.}
		\end{center}
\end{figure}

{\it Interviewer: You were a member of a group who published mathematical textbooks under the pseudonym Nicholas Bourbaki. How did you come to work with Bourbaki?  And what was your role in the group?}

\noindent {\bf Bass:}
I knew that Bourbaki existed, and that 
some people like Serre, Eilenberg, and Grothendieck were members at one time or another.
But mostly they were French. 
So when I was invited, it was quite a shock to me.

My contact with Serre may have had a lot to do with the invitation.
I had nowhere near as broad a kind of mathematical perspective as many of the people in Bourbaki, though my interests are broad.
Perhaps they needed some strength in algebra, especially people who knew something about noncommutative algebra, which I had some work and strength in. This area was not a big tradition in French mathematics. 

I think it was Jean-Louis Verdier who communicated the invitation. He was a close friend and colleague. I ended up being a member from 1970-1982.

Bourbaki is not focused on new mathematics. Instead, it's about a comprehensive and rigorous way of writing a clear and concise exposition of contemporary mathematics. 
This mathematical writing must 
be coherent across a huge landscape, and this leads to some things that are cumbersome from an expository point of view. 

In typical writing, you localize your subject area, and you can adopt conventions that are not needed to be consistent with conventions in another domain of mathematics. But Bourbaki attempts to coherently talk about the whole landscape of mathematics. There are then many cross-referential things that have to be consistent with each other in a way that domain-specific writing would not entail.

Whatever likes or dislikes people have about Bourbaki's style, the writing is done by mathematicians who in their own professional life are gifted expositors. For instance, Serre is a beautiful writer, and very lovely to read, but he also writes things in Bourbaki that fit the Bourbaki style. Writing in Bourbaki is different kind of project. It doesn't have pedagogy so much as logical order and coherence and rigor in mind as the governing feature. It has a conciseness that comes from a severe, minimalist sort of elegance. In some cases, like Lie Theory, I think Bourbaki's exposition is both original and beautiful.

To the extent that I have a kind of strong allegiance to mathematical proof, I guess the style of my writing, even though sometimes careless, sort of fit that style.

One interesting thing about Bourbaki is its working culture. 
At that time, it was a small group, with varying size that involved up to a dozen people or so. 
They had this big architecture of the whole enterprise of mathematics in view. The early books are about foundational structures, for example set theory, algebra, and topology. Then the later chapters treat higher order structures like functional analysis, manifolds and Lie groups where different foundational structures intertwine and are synthesized. There are various chapters that have to be either written or maybe revised and updated, and they all fit into this big architecture.

So how does this writing get accomplished? The nature of the work is largely social. There's a chapter to be written. It's either new or it might be a revision. The group collectively has some conception of what the chapter is supposed to be about, what the overall content is, what it will emphasize. And then some{\it one} -- an individual group member -- is assigned to write a draft of that chapter.

Their funds are modest, they come from publication of these volumes, and they all go into supporting the work. But at least the conditions of the work, the ambience, and things like that, they're not luxurious, but they're very comfortable and pleasant.

\begin{figure}[ht]\begin{center}
	
		\includegraphics[scale=3.00]{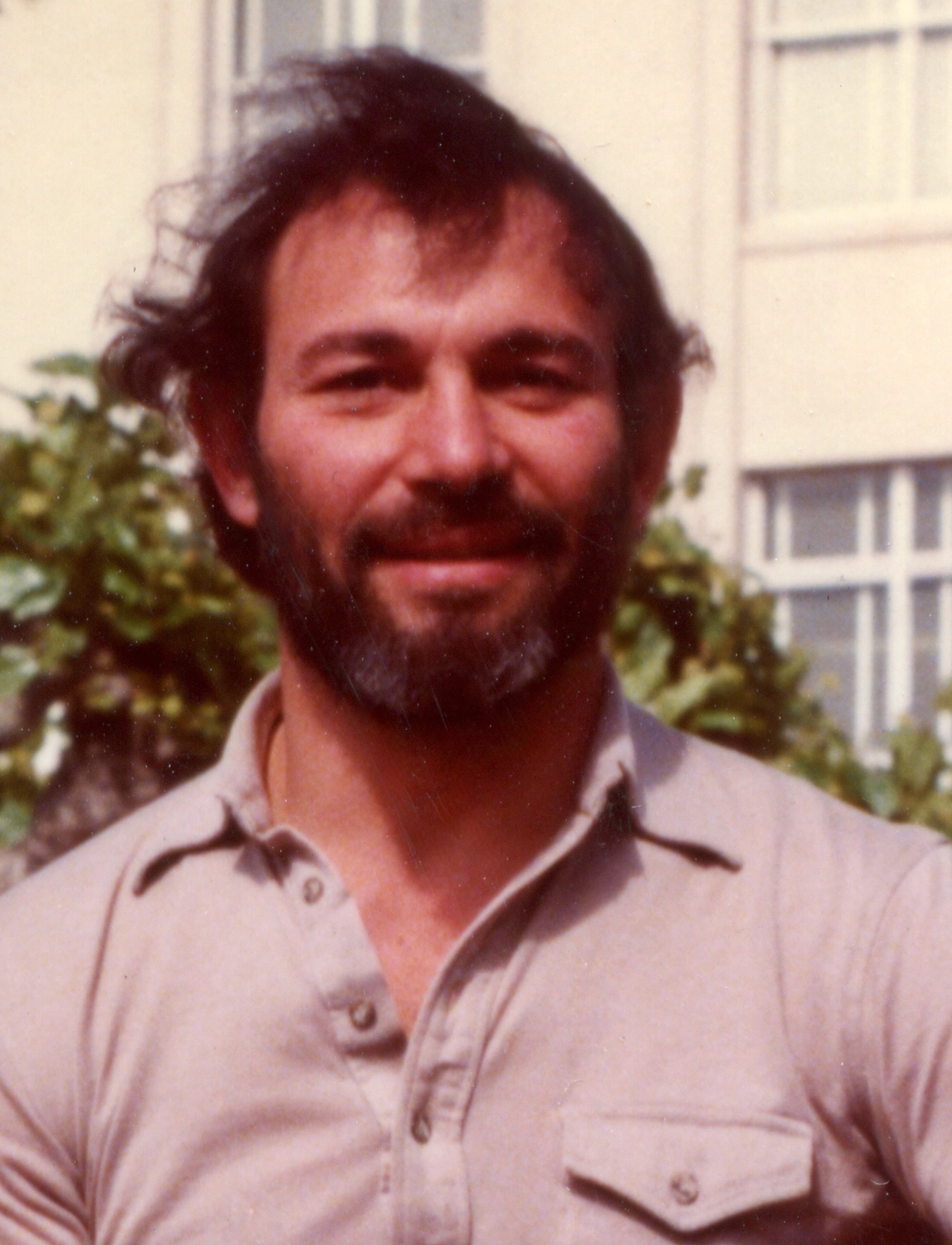}
        \caption{Hyman  visiting UC Berkeley in 1978.}
		\end{center}
\end{figure}

At that time, they met together, in a ``congr\`es'' roughly three times a year. They went to a remote place -- very nice, but quiet, and decidedly not very public. Other people present had no idea who this weird group (of men) are, or what they were doing. The remote locations included Provence, Corsica, and Lake Como, and we stayed in some comfortable auberge and ate very well. For Corsica we took this fast train south from Paris, and had a luxurious meal on board.

The group was very animated but relaxed. Everybody was very casually dressed. They might have been wearing shorts or whatever. We sat around and someone read, line by line, the draft of the current version of a chapter.  The reading traded off at various stages when a reader needed relief. The draft was chewed apart. People raised all kinds of critiques and suggestions – mostly things that they didn't like. The author, of course, was present and participated in this sometimes brutal feedback. The only compensation -- detailed notes are recorded for whatever had to be changed, how things should be done differently -- is that the draft author didn't have to rewrite it.  The next draft was assigned to someone else. A chapter may have gone through many, many cycles only to eventually not die but go {\it au frigidaire} [in the freezer]. In other words, it may be put into long-term storage, possibly forever.  

It's was a very systematic and very exacting kind of process. And, of course, it was anonymous. Other scientists wonder how all this really high level intellectual effort, and time, without personal recognition or compensation, can work\footnote{See `The Bourbaki Gambit' by 
Carl Djerassi, 1996, Penguin Publishing Group}. 
An often publicly discussed question was, ``Why isn't category theory in Bourbaki?''  Well, it was difficult to figure out a way to do that without completely revising Bourbaki's  foundations in set theory. 
Deligne wrote a very nice draft on category theory, and these ideas are {\it au frigidaire}. Interestingly, while Bourbaki does not define functors, he permits himself to use the expression ``proprieties fonctorielles'' in some section headings (though not in the text).

In the discourse of the congr\`es, certain personalities stood out. Jean Dieudonn\'e, Pierre Cartier, and Michel Demazure had pretty strongly developed individual personalities and views about the mathematics. 

My participation in Bourbaki was personally a unique kind of experience, and, by design, largely invisible to the rest of the mathematical world. 

\smallskip
\noindent {\it Interviewer: This does sound like quite a bit of work. So in Bourbaki, the individuals perhaps become famous, but do not receive credit for what they contribute to Bourbaki.}

\noindent {\bf Bass:} What drew me to mathematics was the beauty of mathematical thinking. I'm a very slow thinker. I've always been surrounded by people that were much smarter, much quicker, had a bigger kind of vision of the field than I did. What I valued was to be able to be in the presence of people like that. I took a lot from Bourbaki.
The exposure to that milieu and to those people, was something precious.

For example, I would say that  probably the best piece of work I did was the paper on the congruence subgroup problem. It has 3 authors \cite{BassMilnorSerre1967}. I actually drafted the paper myself, and I can imagine that the other two were not entirely happy with the way it was put together; it could very naturally have been either two or three different papers. For example, the contribution that Serre made came late in the development. I made decisions about some of the terminology which probably would have been better attributed to Milnor's work.

But what I valued about the work was not so much what I did and what I wanted credit for, but the fact that the whole package together was such an impressive and satisfying assemblage of ideas. That to me was most important. I was awe-struck by the fact that this natural question about matrix groups was not only deeply connected to explicit general reciprocity laws in number theory, but that solution of the problem would have required discovery of those laws had they not already been known

For example, I don't think Serre and Milnor, even though they both appreciate each other's work, had ever written anything together. 
So in this sense the ideas that they each contributed are more complementary than intersecting. The whole coherent package was to me  something that I appreciated more; it was not so much about individual credit. It was a question about the aesthetic of the ideas being collected together. It might have been actually more helpful to have separated out parts of it in a way. 

To answer, ``Why do you surrender this time and effort to do something like Bourbaki?'' -- That's just part of the fabric of that whole disposition. I never had any self-image of being a kind of big mathematical creator of ideas and things like that, but I really treasured the idea of being in that environment of really beautiful, high level mathematical thinking. That was most important to me.

\section{1990s-present: Mathematics education}

\begin{figure}[ht]\begin{center}
\includegraphics[scale=0.31]{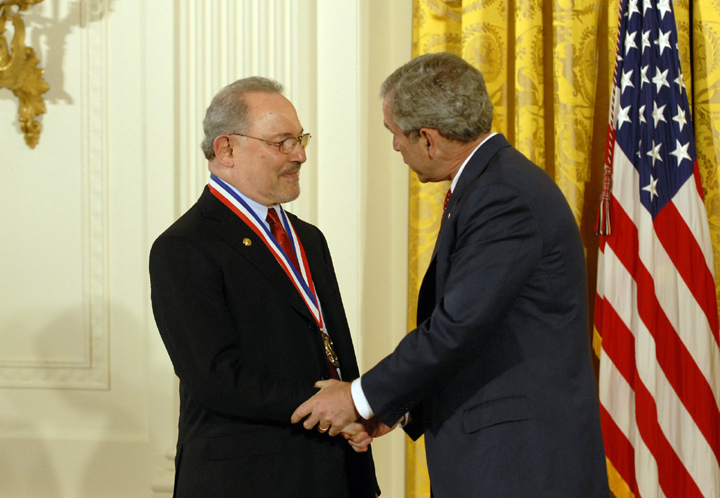}
        \caption{Hyman  receiving the National Medal of Science from President George W. Bush, 2007. Photo by Ryan K. Morris, National Science and Technology Medals Foundation.
        Courtesy: National Science Foundation.
        }
		\end{center}
\end{figure}
    
\subsection*{Inquiry into mathematical work of teaching}

{\it Interviewer: In the 1990s, you started to become more involved in K12 mathematics education. Your longest collaboration in mathematics education is with Deborah Loewenberg Ball. How did the two of you meet?}

\noindent {\bf Bass:} 
There was a workshop organized by Ed Dubinsky at the Educational Development Center, on undergraduate education.
Deborah and I were among those there.
I remember being hugely impressed with her thinking.
She had a kind of worldliness. As I have been saying, I've always been an admirer of really interesting thinking -- mostly mathematical, but also politics, world affairs, civil rights, justice, and people. 
So that's where we met and
we corresponded about these ideas. 

Later, I had more formal contact with her through the Mathematical Sciences Education Board.\footnote{Hyman served on the Mathematical Sciences Education Board, created by the National Research Council, as a member (1991-1993) and Chair (1993-2000).} We invited her to give talks on some of the projects we were discussing intermittently over time. In this way, we had contact that was institutional in nature, but gradually she was telling me more about the research she was doing on teaching. 

As an elementary school teacher, she was always wanting to improve her teaching, but she was troubled by her teaching of mathematics. She felt that there were some understandings of mathematics that
she didn't have, and she needed to find out more. So she took some math courses, which were helpful. But she was uncomfortable importing a disciplinary perspective of mathematics into the classroom without knowing how this perspective was needed for teaching practice. This led to her inquiry on the question, {\it what is the mathematical work of teaching?} But she wasn't convinced that she would see all the relevant mathematical happenings, even if she were able to view a recording of her own teaching. 

It's not often that you see this kind of intellectual humility. People can be very strong at the things they specialize in. But it takes humility to recognize that when you get into some new territory, your ways of knowing your ideas were only for certain kinds of purposes that may {\it not} apply to a different context.  Your knowledge contributes but does not govern.

She then did something that I also like to do 
when want to understand something better. I like to
find the people that I think know and understand the thing most deeply, and 
I go to them and try to
get better understanding. She decided to recruit some mathematicians,
not to tell her what mathematics is appropriate for teaching, but
to look at the actual practice itself. She asked them to watch videos of teaching, guided by the question, {\it what do you see that is mathematically significant and relevant?} At that point, we had substantial intellectual contact, not only about policy, but also about the mathematical work of teaching. 

In these videos, she was teaching third graders. The data we analyzed were amazing.\footnote{The analysis of this corpus of data, the products of Spencer Foundation and National Science Foundation grants, laid empirical groundwork for the notion of mathematical knowledge for teaching.} The videos were accompanied by detailed transcripts. The student notebooks and work were all photocopied. Student interviews and the teacher's journal were also provided. So it's an amazing corpus of data to use to really try to understand what is involved in teaching practice. Deborah has continued in this spirit for a couple weeks each summer at the Elementary Mathematics Laboratory. But this corpus was from teaching from an entire year, the 1989-1990 academic year, when she was teaching third grade. This was remarkable. Nowadays, anyone can do video on their phone. But back then, it was a formidable enterprise to collect this kind of longitudinal documentation at this scale.\footnote{This enterprise, led by Ruth Heaton, Magdalene Lampert, and Deborah Loewenberg Ball, introduced ``hypermedia'' as an approach to understanding teaching practice \cite{LampertHeatonBall1994}.} 

It's surprising, because even though I
tend to be inclined toward abstract ideas, 
this was visceral. When I looked at kids doing things in those classrooms,
it was kind of a
miniature version of seeing the doing of mathematics. It was really exciting. 
I'm not sure
these same kids with a different teacher would have given the same impression.
But somehow
these kids
were thinking about mathematics
and really
invested in it
in a way that I wish more of my undergraduates would be.

The question was to really get inside of this. {\it What is happening mathematically? What is the work of the teaching that sees and discerns what students see and understand?} To the answers to these questions, mathematicians have something to offer. 
But they cannot tell the whole story. There's so much about listening, not just to the kids, but also to colleagues and others -- and moreover attending to the dynamics of instruction and its complexity.  I became engaged with those materials and those questions, and I was strongly attracted to that.

\subsection*{Mathematics as ``making believe''}

{\it Interviewer: How would you describe this kind of teaching and learning?}

\noindent{\bf Bass:} 
 Deborah and I wrote one paper where
we tried to
get across that idea that
we could
map this kind of mathematical classroom culture to what happens in the discipline \cite{BallBass2000MakingBelieve}.

These students, these third graders, were reasoning mathematically. They were
speculating, making conjectures, and arguing. The question we explored was, {\it 
what is the kind of structure of that environment that corresponds to what happens in the discipline?}
In this paper, we say that this structure is ``making believe.''

The phrase ``making believe'' is a common phrase meaning ``pretending'', as in what kids do. But we 
used that as a metaphor for actual
mathematical reasoning, proving, and justification.
Because what do you do when you prove a mathematical theorem?
When is it proven? When is it accepted as part of the knowledge and the
discipline? It's when you have convinced experts in the discipline -- when you have made them {\it believe} it. When the referees say, ``Yeah, this is right, and it's publishable.'' In other words, when you prove a mathematical theorem, it's a process 
of making people believe something.

When you want to prove something, you can say $B$ is true because $A$ implies $B$. What about $A$? You can say $A$ is true because $C$ implies $A$. What about $C$? Where does that stop?
For a mathematician, it stops with axioms, definitions, or prior established results.
You need some things 
that are accepted without need for proof.
What those are varies, of course, depending how expert people are and what the community is. 
We called this the {\it base of public knowledge}. This shared knowledge is not {\it a priori} clear. It must be determined empirically. It is about the common assumptions for that community, such as the definitions you've agreed upon, and the granularity of detail needed. For instance, if you are in high school, if somebody does a routine numerical computation, you don't ask them to justify why that works.
But if it's a sophisticated geometric argument, you want them to explain that in detail. 

In  mathematics we have a notion of base of public knowledge. 
If you're a logician, you're going through the foundations of mathematics. In analysis, you 
don't ask questions about the nature of the real numbers. It's foundational knowledge that you assume. But it's not precisely determined; there has to be some tacit agreement about the starting points.

In other words, each community has its axioms and definitions. And we have theorems that we're allowed to quote that are part of the literature.  

How do you make this base accessible? 
In the third grade classroom, there 
were posters on the wall in the room.
Some posters were about a student contract -- what students agree to do about the norms of discourse and things like that.
But other posters were mathematical. Some of them were challenged later, but at the time they were posted, the students had provisionally reached some consensus that the poster statements were true.
So these posters were the literature. They were a record of what has been achieved.
There were also rules of discourse, how things got challenged, and what the process of argumentation looked like.

To do the work to describe the structure of the classroom that made mathematical work possible, we did not just recount the results and norms. Instead, we gave transcript accounts based on the corpus of data from the 1989-1990 year. Our analysis includes the actual statements of students to each other and what they contested. We identified elements of the process of
mathematical exploration, including the development of mathematical language, and the reasoning of justification.

For instance, in January 1990, the students were working with even and odd numbers. But she (Deborah) didn't start with defining these numbers to start. Instead, she gave the problem:

\begin{quote}
{\it Erasers cost 2 cents, and pencils cost 7 cents. You can buy pencils and erasers. How many different combinations of erasers and pencils can you buy if you want to spend exactly 30 cents?} 
\end{quote}

The point is that the number of pencils has to be even, because they cost 7 cents, and 7 is an odd number. The erasers cost 2 cents each, that 2 is even. So you need an even number of these pencils for the pencils and erasers to add up to exactly 30 cents.\footnote{For an account of this episode, see Ball \& Bass (2003)\cite{BallBass2003Reasonable}.}

And that's how the students discovered there two kinds of numbers, even and odd. 
With those notions, 
they began to explore, what happens when you add two of these kinds of numbers. 
They began to conjecture: 
\begin{align*}
\textnormal{even} + \textnormal{even} & = \textnormal{even} \\
\textnormal{even} + \textnormal{odd} & = \textnormal{odd} \\
\textnormal{odd} + \textnormal{odd} & = \textnormal{even}. 
\end{align*}
Then two students said, we don't agree. One of the students gave an argument that they thought was true. 
But these two students then said, we thought about it, and we don't think you can show that it always works.
Deborah asked them, why is that? And the students said, numbers go on and on forever, and even and odd numbers go on forever, so you can't show that it always works.

It was earth shaking.

These students were right, in the sense that all their proof methods had until this point been done by checking cases.
Earlier in the year, the problems were finite in scope instead of infinite.
They could, by checking all the cases, prove that they had treated all the cases that would come up. 
So these students were wrong that you can't prove it in all cases, and they were right that you couldn't prove it with the method of checking cases.

\begin{figure}[ht]\begin{center}
		\includegraphics[scale=0.49]{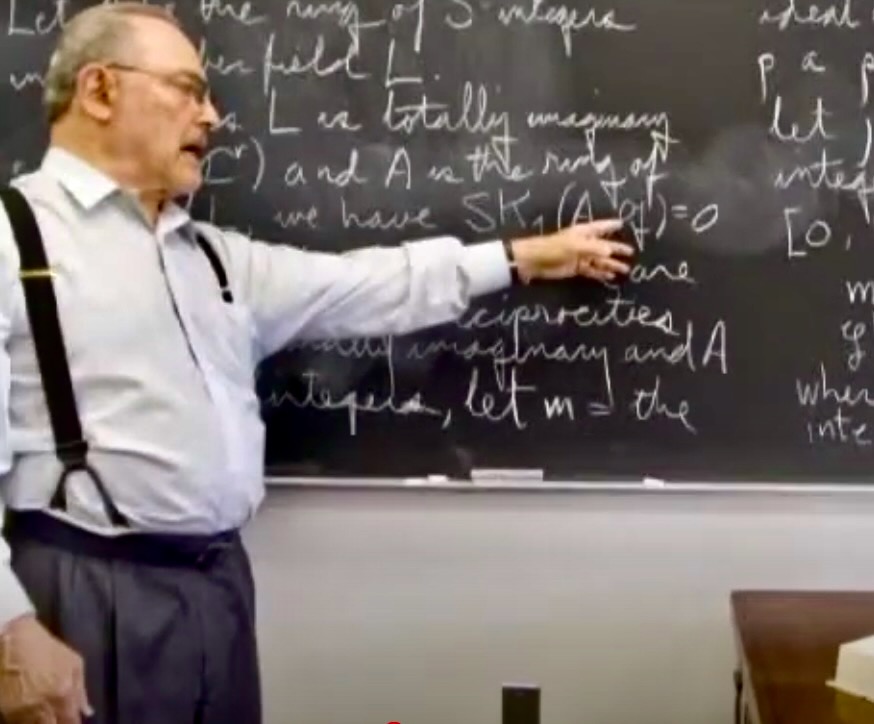}
        \caption{Hyman at the University of Michigan, 2006.}
		\end{center}
\end{figure}

What they were doing was not only
thinking mathematically -- they were actually building up the very infrastructure of proof itself.
They didn't have a formal definition of odd numbers.
But if numbers meant collections of things, and an odd number when you divide it up into pairs, there's always one left over.
And then, of course, when you put two odd numbers together, the two ones leftover can make another pair.
So how did that argument surpass this challenge? It's because 
even though there are infinitely many cases, the definition itself of odd number
is infinite in scope.
It's infinitely quantified. All odd numbers have, when you divide them up into pairs,
one left over.
So the infinite scope of the conclusion 
is possible because it's already
embedded in the very definition.

This is an impressive evolution of mathematical thinking.
And it wasn't programmed. It's not
scripted in some curriculum that tells them to behave this way.
It is a kind of teaching and learning that happens over time. This was January. It's the middle of the school year.

In September, the teacher was still asking them questions like this, and their mathematical thinking it was active, but it was much less evolved than it was by January.
In September, they weren't accustomed to or 
even being confronted with questions like that before. That was new to them. They were able to say interesting things and
begin to develop their curiosity about them. But they had not yet developed the culture that made those September children become the January children that did these amazing things.

That's the importance of having these records
of the instruction through
the entire school year. All the classes were videotaped.
You can see what are the kinds of moves that
supported this development. 
That's not the kind of thing that's easy to capture, even in writing.
It was powerful to be able to 
look at the
primary data.

\subsection*{Building capacity in mathematics teaching}

{\it Interviewer: What conditions would make it possible to learn this kind of teaching?}

\noindent {\bf Bass:} 
 There's something deeply insightful about this kind of teaching. One question is, {\it can you reproduce that kind of teaching
at scale?} How would you build a teacher education program
that would reproduce that kind of thing? You're not going to easily produce another Deborah.
But Deborah herself is very much interested in
creating the capacity
for
those kinds of qualities in teaching.

The approach that she and other like-minded people, such as Pam Grossman, was to look at
how skilled professionals are formed in other areas of work and take lessons from that. How do you become a
surgeon? How do you become
a pilot?
How do you become a plumber?

One idea they have is the
{\it decomposition of practice}.\footnote{Grossman and colleagues (2009)\cite{GrossmanEtAl2009} introduce this concept in a comparative study of training programs in seminaries, clinical psychology, and teacher education.}  When you have a hard, complicated problem, you
can sometimes try to break it up into smaller problems,
which collectively,
maybe with some additional work, can be put together to tackle the whole problem.
In a mathematical proof, we call those lemmas. They're pieces of an argument that you can isolate.

Deborah organized a group in the Michigan faculty, especially when she became Dean. This 
a group that was interested in identifying and developing what they eventually called {\it high-leverage practices} -- things that teachers do daily
as part of their work, that can really impact instruction.

The hope was that these were pieces that teachers can individually become skilled at, and eventually put them together.
It's like if you are a tennis player, you want to know forehand, backhand, and serving, and 
different things like that. But even if you do all those things well, that doesn't mean you'll play well. You've got to coordinate these skills and synthesize them somehow.

\begin{figure}[ht]\begin{center}
		\includegraphics[scale=0.54]{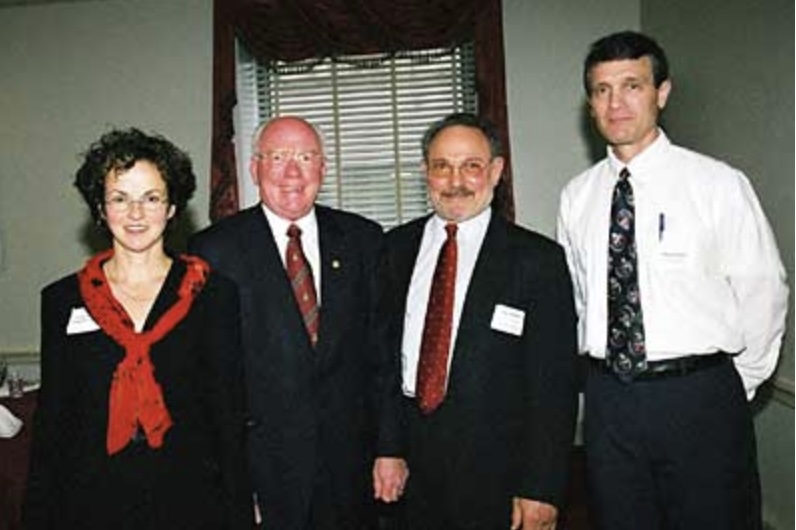}
\caption{ Left-to-right Deborah Loewenberg Ball, Rep. Vernon Ehlers,  Hyman Bass and Roger Howe at the AMS lunch briefing for Members of Congress and their staff on Capitol Hill, July 26, 2001. 
       }
		\end{center}
\end{figure}

Over a period of a couple of years, this group came up with the list of approximately one hundred practices. Now that was unwieldy.
And they were not independent of each other. The next step was to distill these practices into things that were really crucial, that were central to teaching, and that had a big impact. It came down to a list of 19 high-leverage practices. For example, one of these practices was ``leading a mathematical discussion''. This list is continuously modified, and either reduced or expanded. She then set up an institute called TeachingWorks, at Michigan, that does professional development around the country for developing these various
high leverage practices.

Another question in this area is, {\it what does it take to become a teacher when you have to get licensed?}
In the US, licensure is by state. That's not unreasonable in itself. 
For example, in law,
you have different practices in each state, because the laws actually are different.
But in medicine, there are national exams. Doing a heart surgery
in Michigan or in Alabama is not really different. The heart, the body, is generally the same.
Similarly, fractions
don't change when you cross state boundaries.
So the idea of national standards in mathematics
make sense.
But the
political structure in this country is not amenable to
doing things that way. 
At one point Deborah convened a meeting about the idea of a national licensure for teaching mathematics. The idea had a lot of merit, but it's not clear that it would have enough political traction.

\begin{figure}[ht]\begin{center}
		\includegraphics[scale=0.79]{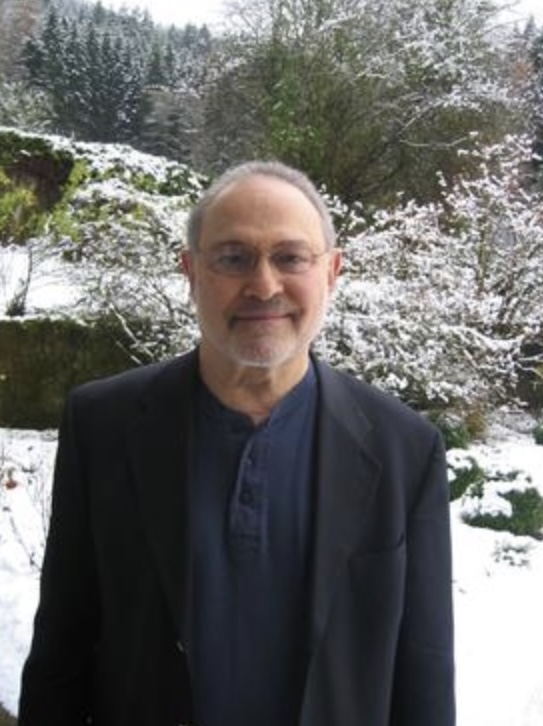}
        \caption{Hyman at a meeting on Professional Development of Mathematics Teachers at Oberwolfach, 2007, photo by Ingeborg Pietzko.}
		\end{center}
\end{figure}

\subsection*{Connection-oriented mathematical thinking}

{\it Interviewer: When did you decide to turn your intellectual focus to math education?}

\noindent {\bf Bass:} 
It was more evolutionary, and age had something to do with it as well. 
Doing serious mathematics is very demanding, not just on time, but also on 
the things you think about when you're
taking a shower or walking.  Memory also plays an important role, and that's harder with age.
It was hard to do serious mathematics at the same time as doing serious thinking about education, both at full strength. 

The work I
did with Deborah
was mainly about the elementary level, but still mathematically sophisticated. 
In the area of elementary mathematics education, much of the work has gone in the direction of thinking about culture, equity, and human behavior. These aspects are important, and absolutely essential. 
But these are also areas where I feel that I have less expertise.

I've never been a quick thinker, even in $K$-theory, where I was heavily involved from the beginning.
Eventually, once the subject became
highly active, there were people that were running with it far faster
than I could, and it took on
a life that was much bigger, with lots of ramifications that were well beyond my expertise.

In education, there are a lot of things I'm interested in, and even have some insights. 
But I know that there are lots of people in the field who've thought about it for longer, and probably with more insight.
I don't feel that I'm in a position to say something that
goes beyond
what they are already able to do well.

So I tried to think, where is something that I feel has been under-treated, and where I could have new insight? Lately it's beginning to coalesce around the idea of mathematical connections, about the unity of mathematics. It is about the coherence, in a large sense, of mathematical ideas that should be part of mathematical education. 

A question here is, {\it if you think about 
mathematics
as a body
of knowledge and skills that people should learn, 
how do we structure an educational system to
equip people with those resources?}
In mathematics, there's analysis and specialization in contrast to synthesis, and there's decomposition and composition.
You can break mathematics down into parts, and decompose it --
that's the analysis.
Doing this, we can identify
coherent chunks of the discipline that you would teach
and give names to them like number theory, algebra, geometry,
probability,
calculus, whatever. We have this core structure which represents a decomposition of the discipline, and it is how our education is structured. 
 But what about synthesis?

I've looked at a lot of curricula and teaching. Even in undergraduate teaching, 
students experience
these courses as islands
that are
disconnected from each other. It's a siloed structure. And this has effects. 

When students get a problem, the first thing they do is typecast the problem.
This is a number theory problem. This is a geometry problem.
This is a calculus problem.

Then, if the problem solution requires that they deploy
resources from more than one domain, 
and if it's not in the type that they categorize the problem,
they unconsciously don't allow themselves to cross that boundary and to another territory.
They haven't been told not to do that. They just don't do that.

It's like when people first learn geometry, no one says, {\it you're not allowed to draw auxiliary lines in that figure}. 
But it is a cultural impediment. So there are lots of signals that suggest that even when students have had a very good
education, and they can show highly developed skills in each area, 
they still don't come away with a sense of the
kind of coherence of mathematical enterprise as a whole. 

Of course, in the outside world, problems are multidisciplinary, and they cross boundaries. Even within mathematics, this is true. 

One message that I think about from cognitive science is that
how much knowledge you have, or what knowledge you have, is not sufficient. What's more important  is how well networked is
the knowledge you have, in other words, is how well developed are the connections between its parts.
There are various pieces of evidence for this, and one from 
cognitive psychology is related to a concept called {\it transfer}.
If you gain knowledge of $A$, can you apply that to $B$,
when $B$ is more or less the same structurally?\footnote{For instance, if you learn negative numbers on a number line at school, can you use this to make sense of money exchanges?}
Transfer is about whether you can export knowledge to structurally similar problems in different contexts. 
Well, experiments show that transfer very rarely occurs.

So I've been thinking about, 
{\it How can you teach people to make mathematical connections?} 

Of course, there are many kinds of connections. Some of them are fairly obvious. If you give a drill sheet of exercises, those problems are all related, and students see that they're related. 

\begin{figure}[ht]\begin{center}
		\includegraphics[scale=0.22]{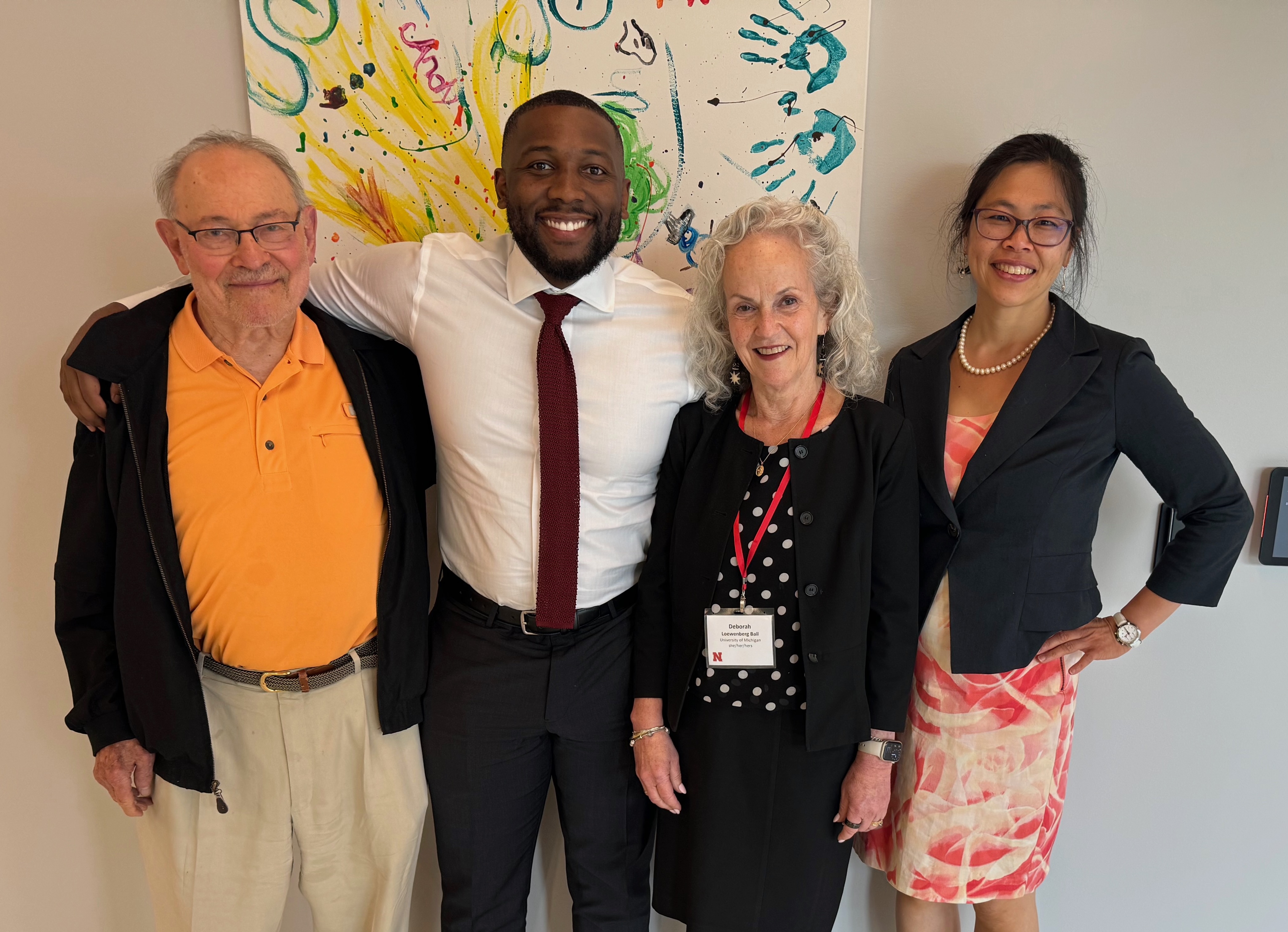}
        \caption{Left-to-right Hyman Bass, Charles E. Wilkes II, Deborah Loewenberg Ball and Yvonne Lai at the University of Nebraska-Lincoln, May 2024.
       }
		\end{center}
\end{figure}

I'm interested in connections that are not obvious. 
Connections that have some depth and subtlety.
One reason that I'm drawn to this is that 
some of the big breakthroughs in math, and also in science, are when people realize that two things that seem unrelated, that perhaps are in
different parts of the subject,
are actually
connected in some deep way.

For example, Poincar\'e's realization that 
the transformations of Fuchsian functions he was looking at were identical with the transformations of non-Euclidean geometry came to him
as he was boarding a bus.\footnote{``…just as I put my foot on the step, the idea came to me. Nothing in my former thoughts seemed to have prepared me for it: the transformations I had used to define Fuchsian functions were identical with those of non-Euclidian geometry. I had no time to verify the thought, because the conversation started again as soon as I sat down, but I felt absolute certainty at once.'' -- Poincar\'e, 1908, p.~388 \cite{Poincare1908}} 
Or take the Langlands program. 
It relates number theory and 
Galois representations to harmonic analysis. It's not a theory; it's a program that people will spend another several decades on.

Mathematics is about finding patterns. Mathematicians find in their own practice
that they're seeing a pattern when it seems like, in diverse contexts, they're doing essentially the same work
over and over again. 
When that happens, 
they ask, {\it what is
the `thing' of which these are all special cases?}
That's a kind of conceptual question, because it's a new concept
that doesn't yet exist.

Abstraction is not such an exotic thing. Numbers are abstractions. The number 5. Where's the number 5? Can you point to it? No.
But you can hold up 5 fingers, and that's an example of 5 things. And I can say, if another thing is an example of 5 things, I can concretely implement that I can show that there's a 1-to-1 correspondence between the elements. Saying the two things are the same is concrete.
But what is
giving a name to that sameness? This is creating something new, a new concept that didn't exist before.
And that's generally what abstractions are about.

My thought was that when you build a theory, it's creating a conceptual frame
that encompasses lots of things that seem different as being cases of the same thing.
I've tried to approach this idea from a lot of different directions. 

Now, we can't teach
students to have the kind of visionary insight that brought the Langlands program or Poincar\'e's insight, because, first of all, you have to know a lot of different things to do that.
But 
what I've been trying to develop lately is what is called 
``connection-oriented mathematical thinking".

And I can give lots of accessible examples of that.\footnote{Some examples may be found in \cite{Bass2017}, for instance when Hyman  demonstrates a general measure theoretic structure that underlies a set of problems on mixing liquids, discrete counting, and area.}
The point is that discerning and using connections is
not something we can expect without instruction. So the question becomes, {\it 
how can you teach connection-oriented mathematical thinking?}
Those are the ideas I've been exploring lately. I have a sabbatical this coming semester.
I'm hoping to be able to write something about it.

In the transfer experiments in cognitive psychology, participants would be given a pair of so-called isomorphic problems, say 
Problem $A$ and Problem $B$.
The experimenters would teach participants Problem $A$ very well, then give them Problem $B$.
The experimenters wanted to know whether participants could recognize the structural equivalence and mobilize knowledge of $A$ to solve $B$.
They find this rarely happens. 

My focus is {\it not} on transfer. Instead, I foreground {\it the making of the connections}. The point is that even if you may be able to solve $A$ and $B$ independently, you still may not see the connection between them. And, sometimes, even when problems seem like they're equivalent, they aren't. 

Even with fairly elementary problems, it can be challenging to sort or type them according to how 
structurally connected they are, in particular whether they are actually equivalent or isomorphic.
There are various materials that I have been developing to design instructional activities for this kind of thinking.

\section{Relationships among mathematicians and educators}

\begin{figure}[ht]
\begin{center}		\includegraphics[scale=0.15]{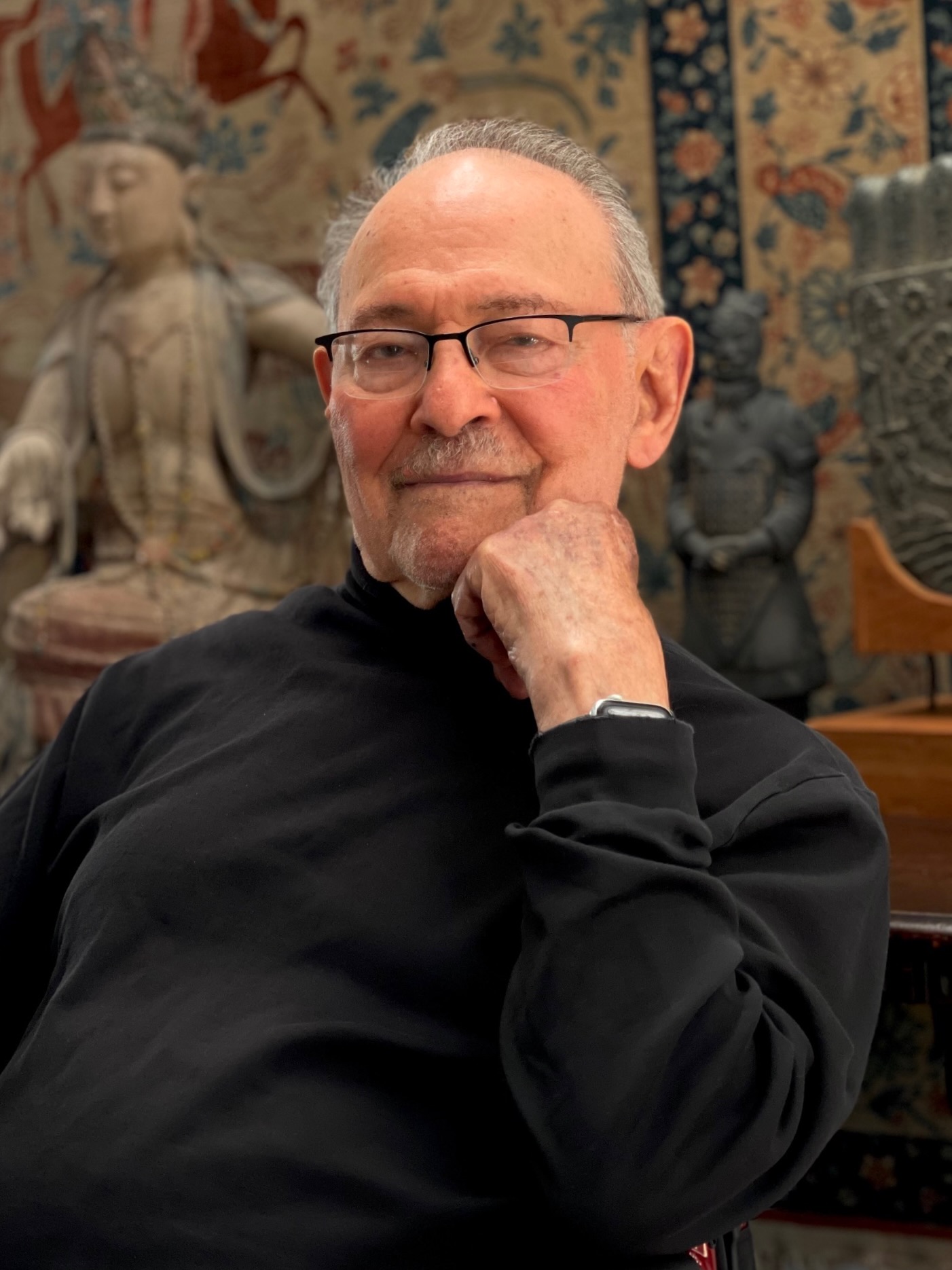}
        \caption{Hyman at home in Dexter, Michigan, 2024. Photo by Gabriella Bass.}
		\end{center}
\end{figure}
    
{\it 
Interviewer: In the Math Wars of the 1990s to 2000s, you served in various governance positions, including in the AMS and the National Academies.  What do you think was the largest impact of these Math Wars?}

\noindent {\bf Bass:} 
I think the Math Wars were
somewhat artificial. Still, they set things back for a long time.

It's a complicated story. One of the things that 
mathematics and mathematicians suffer from
is a kind of
elitism, or even arrogance. 
In the Math Wars,
some mathematicians have tried to 
assert their authority because they proved some theorems, and they have a position in a good university. Therefore
they're not subject to question.
It's like society decides that if you have some really highly developed skill and some technical area, that therefore you're wise and have insights and things to offer to almost any area of human concern. 
Even the assumption that when engineers are out of employment, because some industry has shifted, that they are qualified to be teachers, is part of this. 

Mathematicians are very particular and
fastidious about
treating the subject matter with rigor and integrity. 
But, interestingly enough,
even
good mathematicians don't always agree about the best way of
rendering these ideas -- 
much less knowing about the entailments of teaching.

So
the fact is that there is often a general kind of contempt for teachers.
People may not express this directly. It's just that teachers have no standing. There is a stance among some mathematicians that teachers are not to be listened to, but rather that teachers should listen to us. 
That's a real barrier to progress.

How are you going to elevate the level and standing of the teaching profession? This is complicated.
There are lots of problems even within education to be able to make that happen.
For example, developing high level practice involves high level critique.
If you're an artist or an architect, there's a practice of `crit' by peers. In  studying medicine, there are rounds. 
In an art studio, other artists gather around your draft painting, or something you’ve done, and they ask really good, insightful questions. Why did you do this? 
What do you think about
positioning something differently? What about the color, composition, \dots ? 
Artists and architects and surgeons are accustomed to that. They're serious practitioners who are listening to peers.
There are many fields where that's done very effectively. 

But teaching is very private. 
Teachers often feel unsafe when other people come in and
try to evaluate what they're doing -- and for good reason. The way evaluation has been conducted has been so punitive and unconstructive, and has sometimes been done by people who lack adequate understanding of the work. 
So there's much that has to change.
Mathematicians have not helped make it
safe for people to develop and improve
culture.
If you want to improve the practice of a huge profession, you're not going to advance
the cause by deriding the members of that profession. 

Of course, debate is appropriate. Education has to evolve, and 
there are good ideas in the field such as the concept of continuous improvement and looking at things at the system level (e.g., \cite{ParkHironakaCarverNordstrum2013}).

The work that Deborah does is very much on the ground at the level of practice.
Nothing
can be
fully consummated without
seeing what practitioners do
with it at the ground level.
But there also has to be
thinking more at the systemic level.  But, ironically, the actual work of teaching is the most understudied part of education.

\smallskip

\noindent {\it Interviewer: What is your wish for the future of interactions between mathematicians and educators?}

\noindent {\bf Bass:} 
Well, overall, I would say harmony.
I think that
a lot of progress has actually been made on this. 

For mathematicians to take education more seriously -- not just respect for K12, but even at the university level. 

Not everybody is expected to do that
with equal intensity.
But much more respect
and standing has to be given to the people in their departments who do that work. Mathematicians need to recognize this work as
equally important to, and even synergistic with,
the research that they do, and, in fact, as something that will enrich their research.

\bibliography{BassLong.bib}

\end{document}